\documentclass[11pt]{amsart}
\usepackage{amssymb}
\usepackage{amsmath} 
\usepackage{amsthm} 
\usepackage{calc} 
\usepackage{verbatim} 
\usepackage[dvips]{graphics}
\usepackage[dvips]{epsfig}
\newtheorem{thm}{Theorem}[section]
\newtheorem{lem}[thm]{Lemma}

\newtheorem{cor}[thm]{Corollary}

\newtheorem{prop}[thm]{Proposition}

\renewcommand{\qed}{\hfill$\Box$\medskip}

\newcommand{\ignore}[1]{}

\newcommand{\Area}{\mathrm{Area}\,}

\newcommand{\interior}{\mathrm{int}\,}

\newcommand{\dist}{\mathrm{dist}}

\newcommand{\eqdef}{\stackrel{\mathrm{def}}{=}}

\newcommand{\id}{\mathrm{id}}

\newcommand{\norm}[2][{}]{\left\|#2\right\|_{#1}}

\newcommand{\pair}[2]{\left\langle #1,#2 \right\rangle}

\newcommand{\R}{\mathbf{R}}

\newcommand{\C}{\mathbf{C}}
\newcommand{\cp}{\mathbf{P}}
\newcommand{\Z}{\mathbf{Z}}
\newcommand{\N}{\mathbf{N}}

\newcommand{\ti}{\tilde}

\input epsf.sty
 
\begin{document}
\newcommand{\gr}{\lambda}
\newcommand{\pp}{\cp^1\times\cp^1}\
\newcommand{\cB}{\mathcal{B}}
\newcommand{\cC}{\mathcal{C}}
\newcommand{\cD}{\mathcal{D}}
\newcommand{\cF}{\mathcal{F}}
\newcommand{\cI}{\mathcal{I}}
\newcommand{\cK}{\mathcal{K}}
\newcommand{\cO}{\mathcal{O}}
\newcommand{\cP}{\mathcal{P}}
\newcommand{\cS}{\mathcal{S}}
\newcommand{\cT}{\mathcal{T}}
\newcommand{\cW}{\mathcal{W}}
\newcommand{\bc}{\mathbf{c}}
\newcommand{\bb}{\mathbf{b}}
\newcommand{\Per}{\mathrm{Per}}
\newcommand{\SPer}{\mathrm{SPer}}
\newcommand{\total}{\dot}
\renewcommand{\eqdef}{:=}

\setlength{\textwidth}{6.5in} 
\setlength{\textheight}{9in} 
\setlength{\headheight}{0in} 
\setlength{\evensidemargin}{0in}
\setlength{\oddsidemargin}{0in}
\setlength{\topmargin}{0in}

\title{Real and Complex Dynamics of a Family of \\ 
        Birational maps of the Plane:\\
 The Golden Mean Subshift}


\author{Eric Bedford \& Jeffrey Diller}
\address{Department of Mathematics\\
Indiana University\\
Bloomington, IN 47405}
\email{bedford@indiana.edu}
\address{Department of Mathematics\\
         University of Notre Dame\\
        Notre Dame, IN 46556}
\email{diller.1@nd.edu}

\subjclass{}
\keywords{}


\maketitle
\markboth{}{}

\section*{}
{\footnotesize
{\bf Contents}

0. Introduction \dotfill 1

1. Filtration \dotfill\ 4

2. Coding and Rectangles \dotfill 8

3.  Complexification and Intersection Theory \dotfill 12

4.  Structure of Rectangles \dotfill 14

5.  Invariant Conefields; Boundaries of
Rectangles \dotfill 20

6.  Periodic Points \dotfill 23

7.  Uniform Arcs and One-sided Words \dotfill 26

8.  Conjugacy with the Subshift \dotfill 29

9.  Parabolic Basin; Nonwandering Set \dotfill 33

10.  Stable Manifolds and Laminar currents 
\dotfill 35

11.  Parameter Values $-1<a<0$ \dotfill 40

12.  The Purely Complex Point of View  \dotfill 42

}
\bigskip

\section*{0.  Introduction}
\label{intro}
\noindent  Much attention has been given to birational mappings which arise in connection with integrable
systems; two recent surveys of this subject are given in [BTR] and [GNR].  For these birational mappings, it is 
of interest to know the behavior of the iterates $f^n=f\circ\cdots\circ f$ as
$n$ increases.  A family which comes from the study of lattice statistical mechanics (see [BM]) is
$$f_a(x,y)=(y{x+a\over x-1},x+a-1)$$
for fixed $a\in{\bf R}$, which defines a birational mapping of the plane ${\bf R}^2$.  (Throughout this paper we
exclude the case $a=-1$ since $f_{-1}$ is affine and thus dynamically trivial.)  Typical of mappings that
arise this way,
$f_a$ is area-preserving in the sense that it preserves a meromorphic 2-form, and is reversible,
which means that $f_a$ is conjugate to $f^{-1}_a$ via an involution.  This family was investigated in a series of
papers by Abarenkova et al.\  [Ab1--5], which describe several numerical phenomena and raise a number of
interesting questions. 

The goal of this paper is to give a precise description of the dynamics of $f_a$ for $a<0$.  One of our
motivations here is to give a first example of pointwise dynamics of a
birational mapping which is, in an essential way, not a homeomorphism of its nonwandering
set.   One of the properties of a rational map is that it can have points of
indeterminacy: at such points the mapping cannot be defined to take a single value, and
the map is said to ``blow up'' such points, assigning whole curves to them. 

One way of dealing with such behavior is to consider the closed relation induced by $f_a$ and to work within the
general framework of topological dynamics (see Akin [Ak]).  This framework, however, does not reflect the rich
structure of a birational map, so we also work with the complexification $\tilde f_a$,
which is a birational map of ${\bf C}^2$.  This allows us to use the tools of
complex analysis, complex geometry, and complex potential theory.  Despite pointwise difficulties,
$f_a$ induces a well-defined map $f_a^*$ on the set of positive closed (1,1)-currents (see Sibony [S]
and Guedj [Gu1]).   In [DF] it was shown that for all values
$a\in{\bf C}-\{-1,0,{1\over 3}, {1\over 2}, 1\}$
 there is an invariant (1,1)-current $T^+_a$ such that $\tilde f^*_aT^+_a=\rho_a T^+_a$, with $\rho_a>1$. 
It follows that the degrees of $f_a^n$ grow
exponentially like $\rho_a^n$, thus confirming a conjecture of [Ab1--5].  Except for a countable set
of exceptional values of
$a$, $\rho_a$ is equal to the golden mean $\phi:=(1+\sqrt5)/2$.

Measures (which correspond  to 0-currents) transform differently from 
(1,1)-currents and thus reflect more closely the pointwise behavior  of the map.  A construction of invariant
measures for a rather  general class of birational mappings, which includes the family 
$\{\tilde f_a\}$ is given in [BD].  In this paper 
we are able to go farther and give a pointwise description of the 
dynamics.  In particular, we describe the behavior of $f_a$ on the 
indeterminacy and critical sets, both of which are ``invisible'' from 
the measure-theoretic point of view.

Since the coordinate functions of $f=f_a$ are rational, we extend $f_a$ to the compactification
$\overline{{\bf R}^2}:=({\bf R}\cup\{\infty\})\times({\bf R}\cup\{\infty\})$.  The point
$(\infty,\infty)$ is a parabolic fixed point for
$f$.  The forward/backward basin $\cB_\pm$ is the set of points where $f^n$ converges locally uniformly to
$(\infty,\infty)$ as
$n\to\pm\infty$.  We prove in Section 9 that the nonwandering set is the complement of the
parabolic basin $\cB_+\cup\cB_-$.

Our abstract model for the dynamics of $f$ on the nonwandering set will be the so-called golden mean
subshift
$(\sigma,\Sigma)$.  That is, 
$\sigma$ is the shift map, and $\Sigma$ is the topological space of bi-infinite sequences of 0's and 1's such that
`1' is always followed by `0'.  The entropy of this subshift is the logarithm of the golden mean $\phi$.  We
connect
$f_a$ with our model system by giving a (multiple-valued) equivariant correspondence $R$ taking
$(\sigma,\Sigma)$ to
$(f,\Omega)$.  Note that $R$ cannot be a topological conjugacy because $(\sigma,\Sigma)$
is a homeomorphism while $(f,\Omega)$ is not. 

We identify rectangles $R_0$ and $R_1$ which cover $\Omega$ and which serve
as a Markov partition.  Any $f$-orbit $(p_n)_{n\in{\bf Z}}$ which lies in $\Omega-(\infty,\infty)$ can be assigned
a coding
$w=(w_n)_{n\in{\bf Z}}$, where each symbol
$w_n\in\{0,1\}$ is chosen so that $p_n\in R_{w_n}$.  By the mapping properties of $R_0$ and $R_1$, it follows
that `11' cannot occur, and thus $w\in\Sigma$.    We adopt the convention of coding the fixed point
$(\infty,\infty)\in\Omega$ by the 2-cycle $\overline{01}\leftrightarrow\overline{10}$ in $\Sigma$.  As it turns
out, we will find it more convenient to work with the set-theoretic ``inverse'' of the coding map.  Let $R(w)$ be
the ``rectangle'' of points coded by $w$.   More precisely, we set
$$R(w)=\bigcap_{n\in{\bf Z}}f^{-n}R_{w_n}-(\infty,\infty),$$
for $w\in\Sigma-\{\overline{01},\overline{10}\}$, and we set
$R(\overline{01})=R(\overline{10})=(\infty,\infty)$.
$R$ defines a semi-conjugacy from $(\sigma,\Sigma)$ to $(f,\Omega)$ in the following sense:  If $R(w)$
does not contain the point of indeterminacy $(-a,\infty)$, then
$$fR(w)=R(\sigma w).$$
There is a unique word $w_*$ such that  $R(w_*)=E$ is a nontrivial interval containing a point of indeterminacy
(see Figure 6).    One of our principal results is
Theorem \ref{finiteomega}, which says that if
$w\ne\sigma^nw_*$, then  $R(w)$ is a single point.  

Now we may further describe the dynamics of $f$: 
$f^2$ acts by translation on the two lines at infinity, so the behavior is
decidedly non-hyperbolic everywhere on $\overline{{\bf R}^2}-{\bf R}^2$ .  On the other hand, the behavior of
$f$ on $\Omega\cap{\bf R}^2$ has many of the properties of an Axiom A diffeomorphism:
\begin{enumerate}   
\item  There are invariant cone fields for
$f$ at all points of $\Omega\cap{\bf R}^2$.  
\item All periodic points, except for
$(\infty,\infty)$, belong to ${\bf R}^2$ and are saddle points.  
\item The  saddle
points are a  dense subset of
$\Omega\cap{\bf R}^2$.  
\item   $f$ is topologically expansive on
$\Omega\cap{\bf R}^2$.  
\item There are stable and unstable
manifolds through every point of
$\Omega\cap{\bf R}^2$; the corresponding  laminations, $\cW^s$ and $\cW^u$, intersect
transversally.   
\end{enumerate}
Let us reiterate that the cone field and expansivity mentioned above are defined only on $\Omega\cap{\bf
R}^2$, which is not an invariant set.  Note, too, that distinct stable manifolds intersect at points of the
countable set
$I_+=\bigcup_{{n\ge0}}f^{-n}I(f)$; see Figures \ref{omegapprox} and \ref{omegapprox2}.

Finally, we draw a parallel with a result of Ruelle and Sullivan [RS] for Axiom A surface mappings; we show
that
${\cW}^s$ and
${\cW}^u$ can be used to construct a ``stable current''
$\mu^+_{\bf R}$ and an ``unstable current''
$\mu^-_{\bf R}$ whose intersection product gives the unique measure of maximal entropy.  These
currents should give a connection between the real dynamics of $f_a$ and the complex dynamics of
$\tilde f_a$ because the invariant current
$T^+_a$ (which has real dimension 2) appears to be the ``complexification'' of the 1-dimensional current
$\mu^+_{\bf R}$.

Let us outline our mathematical approach.  We work simultaneously with the real map $f_a$ and its
complexification $\tilde f_a$.  We
consider the forward/backward iterates of complex lines in ${\bf C}^2$.  Let $\tilde L$ and $\tilde M$ denote the
complexifications of real lines $L$ and $M$.  By the intersection theory of complex subvarieties, we know that
the intersection number of $\tilde f^{-n}\tilde L$ and $\tilde f^m\tilde M$ is determined by the homology classes
of these two sets.  Considering the purely real behavior, we develop a geometric/combinatorial argument which
gives a lower bound on the number of (real) intersection points of 
$f^{-n}L\cap f^mM$.  This allows us to conclude that $R(w)$ is nonempty. This lower bound coincides with the
upper bound given by complex intersection theory; hence all intersection points are real and have
multiplicity one.  The property of having multiplicity one leads to transversality and the existence of cone
fields.  

We believe that the maps $\{f_a: a<0\}$ represent an important subfamily within the whole family $\{f_a\}$.
This may be seen by analogy with the H\'enon family $\{h_{a,b}:b\ne0\}$
$$h_{a,b}(x,y) = (a-x^2-by,x)$$
of quadratic diffeomorphisms of ${\bf R}^2$.  For $a\ll-1$, the map $h_{a,b}$ has
dynamics which are completely transient: all orbits tend to $\infty$.  On the other hand, it was shown in [DN]
that $h_{a,b}$ generates a horseshoe when $a\gg1$.  [HO], using complex methods, were able to obtain a much
larger family of horseshoes.   By Friedland and Milnor [FM] it is known that the entropy of
$h_{a,b}$ is no greater than $\log2$, and thus the horseshoe mappings in $\{h_{a,b}\}$  represent elements of
maximal entropy.  The transition of behaviors of
$h_{a,b}$ as
$a$ passes from
$a\ll-1$ to
$a\gg1$ is seen as illustrating a mechanism for the transition to chaos, with the horseshoe mappings
exhibiting fully developed chaos.  The central position of the horseshoe in this picture is given from the point
of view of the ``Pruning Front Conjecture'' by de Carvalho and Hall in [dCH].  The point of our analogy here is that
the maps
$\{f_a,a<0\}$ have maximal entropy within the family $\{f_a:a\ne-1\}$, and we expect them to play a
fundamental role within this family, as the horseshoes play within the H\'enon family.  On the technical level, too,
we have borrowed from the analogy with the H\'enon family.  It was shown in [BLS] that if $h_{a,b}$ is of
maximal entropy, then the nonwandering set for the complexification is contained in ${\bf R}^2$.  In [BS1,2], it
was shown that a mapping $h_{a,b}$ with maximal entropy may be studied by working with its complexification;
that approach has motivated some of the work of the current paper.

\section{Filtration}
\label{filtration}
\noindent Throughout this paper we consider $f_a$ only for $a<0$, $a\neq -1$.  In fact, we
will assume  that $a<-1$ until we reach Section \ref{smallnega}, where we indicate the
modifications needed to treat the case $-1<a<0$.  
We write $f=f_a$; the inverse of $f$ is given by
$$f^{-1}(x,y)=(y+1-a, x{y-a\over y+1}).$$
The involution $(x,y)\mapsto(-y,-x)$ conjugates $f$ to $f^{-1}$.  The indeterminacy set consists of the points where $f$ takes the form ${0\over0}$ or
$\infty\cdot 0$ and is given by 
$I(f)=\{(1,0),(-a,\infty)\}$.  The critical set is $\cC(f)=\{x=1\}\cup\{x=-a\}$, which contains
$I(f)$.  The critical set for the inverse is $\cC(f^{-1})=\{y=-1\}\cup\{y=a\}$, and the indeterminacy
locus is $I(f^{-1})=\{(0,-1), (\infty, a)\}$, as is seen by applying the involution to
$\cC(f)$ and $I(f)$. 
$f$ is smooth on
$\overline{{\bf R}^2}-I(f)$, and
$f:\overline{{\bf R}^2}-\cC(f)\to\overline{{\bf
R}^2}-\cC(f^{-1})$ is a diffeomorphism.

A calculation shows that $f$ preserves the meromorphic two form
$$
\zeta = \frac{dx\wedge dy}{y-x+1}.
$$
This form has no zeroes; it has simple poles along the lines
$\{x=\infty\}$, $\{y=\infty\}$, and $\{y-x+1=0\}$.  The
union of these lines is invariant under $f$.  One checks directly
that $f$ maps $\{y-x+1=0\}$ onto itself by $(t,t-1) \mapsto (t+a,t+a-1)$,
and that $f$ interchanges $\{x=\infty\}$ and $\{y=\infty\}$ according
to $(\infty,y)\mapsto (y,\infty)\mapsto (\infty,y+a-1)$.  In particular
$f^2$ restricts to translations on these three lines.  It follows that
$Df^2_{(\infty,\infty)} = \id$, so $(\infty,\infty)$ is a parabolic fixed point for $f$.  This
point plays a central role in the dynamics of $f$, so 
we record some information about the behavior of $f$ nearby.

\begin{prop}
\label{approx}
There exists a constant $C>0$ such that $|x|,|y|>C$ implies 
$$
\left(\begin{matrix}
x \\ y
\end{matrix}\right) 
\overset{f^2}\mapsto 
\left(\begin{matrix}
x (1+\frac{a-1}{x} + \frac{a+1}{y} + O(|x|^{-2} + |y|^{-2})) \\ 
y (1+\frac{a-1}{y} + \frac{a+1}{x} + O(|x|^{-2}))
\end{matrix}\right). 
$$
In particular, if we set $m(x,y) = y/x$ (which is negative in $R_0\cup R_1$), 
then
$$
m\circ f^2(x,y) = m(x,y)\left(1 - \frac{2}{y} + \frac{2}{x} 
	        + O(|x|^{-2}+|y|^{-2})\right).
$$
\end{prop}

\noindent The proof is a straightforward computation that we leave
to the reader.

We say that the parameter $a\in\C$ is \emph{exceptional} if
$I(f)\cap f^{n}(I(f^{-1}))\neq\emptyset$ for some $n\geq 0$.  Since
both $I(f)$ and $I(f^{-1})$ are contained in the three invariant lines,
the exceptional values of $a$ are those for which
$(-a,\infty)=f^n(\infty,a)$ or $(1,0)=f^n(0,-1)$.  This happens when 
$a=\frac{n-1}{n+1}$
or $a=\frac1n$
for some integer $n\geq 1$.  Thus no value $a<0$ is exceptional.  If $a$ is not exceptional, then
$I(f^n)=I(f)\cup\dots\cup f^{-n+1}I(f)$.  We use the notation
$$I_+=\bigcup_{n\ge0}f^{-n}I(f)=\bigcup_{n\ge1}I(f^n),\ \ \ \
I_-=\bigcup_{n\ge0}f^nI(f^{-1})=\bigcup_{n\ge1}I(f^{-n}),$$ so that $\overline{{\bf
R}^2}-I_+$ (resp.\
$\overline{{\bf R}^2}-I_-$) is the set of points where the pointwise forward (resp.\
backward) dynamics is uniquely defined.  Thus $\overline{{\bf R}^2}-(I_+\cup I_-)$
is the set of points which are contained in  only one bi-infinite orbit.  The largest
invariant subset of $\overline{{\bf R}^2}-(I_+\cup I_-)$ is
$$\cD_f:=\bigcap_{n\in{\bf Z}}f^n(\overline{{\bf R}^2}-(I_+\cup I_-)) ={{\bf
R}^2}-\bigcup_{n\ge0}\left( f^n\cC(f^{-1})\cup f^{-n}\cC(f)\right).$$ 
Thus $\cD_f$ is the largest set which
$f$ maps homeomorphically to itself.  $\cD_f$ is clearly dense in $\overline{{\bf R}^2}$.

\begin{figure}
\centerline{\epsfxsize4.6in\epsfbox{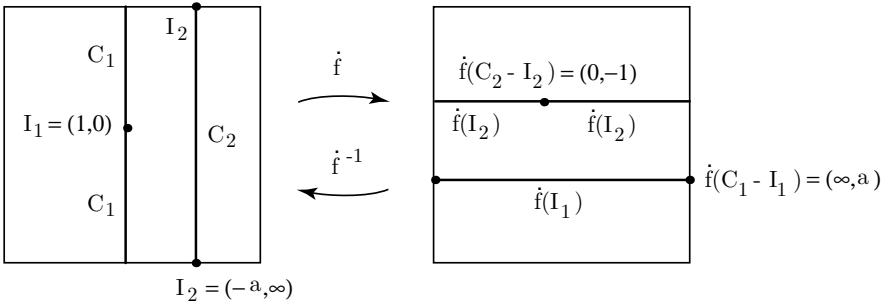}}
\caption{Action of $\dot f$ and $\dot f^{-1}$ on critical sets and points of indeterminacy.}
\end{figure}

We let $\dot f$ denote the closed relation on
$\overline{{\bf R}^2}$ which is obtained by taking the closure of the graph of $f$
restricted to $\overline{{\bf R}^2}-I(f)$.  (We mention Akin [Ak] as a reference for
basic material about closed relations.)  In other words, 
$\dot f$ is the set-valued mapping given by
$\dot f(p)=f(p)$ for
$p\in\overline{{\bf R}^2}-I(f)$, $\dot f(1,0)=\{y=a\}$, and $\dot f(-a,\infty)=\{y=-1\}$.  Figure 1
shows how $\dot f$ acts on $\overline{{\bf R}^2}$; $\cC(f)-I(f)$ is taken to $I(f^{-1})$, and
$I(f)$ is taken to
$\cC(f^{-1})$.  Since $a$ is not exceptional, the operation of passing to the corresponding closed
relation respects the dynamics.  That is, $(\dot f)^n=\dot g$, where $g=f^n$.  Let
$\dot f^{-1}$ denote the closed relation obtained from the restriction of $f^{-1}$ to $\overline{{\bf
R}^2}-I(f^{-1})$.  With
$\cC(f)=C_1\cup C_2$ as in Figure 1, we see that if $q\in C_j$, then $\dot f^{-1}\dot fq=C_j$, for
$j=1,2$.

Let $\cK$ denote the set of compact subsets of $\overline{{\bf R}^2}$.  The relation $\dot f$ induces
a map $\dot f:\cK\to\cK$.  Since $\dot f$ is a closed relation, it is upper semicontinuous on $\cK$.
That is, if closed sets $S_j$ decrease to $S$, then $\dot fS_j$ decrease to $\dot fS$.  

There is a second induced map
$f:\cK\to\cK$, where $f(S)$ is defined as the closure of $f(S-I(f))$.  This map is neither
upper or lower semicontinuous.  For
$S\in\cK$ we have
$$f(S)\subset \dot f(S)\subset f(S)\cup\cC(f^{-1}).$$
Again, since
$a$ is not exceptional, the $n$th iterate of $f$ as a map of $\cK$ coincides with the mapping
of
$\cK$ induced by
$f^n$.  We also have an induced mapping
$f^{-1}:\cK\to\cK$. If
$S$ is a compact set which is the closure of its interior, then
$f(S)$ is also the closure of its interior, and $f^{-1}f(S)=S$.

\begin{figure}
\centerline{\epsfxsize4.6in\epsfbox{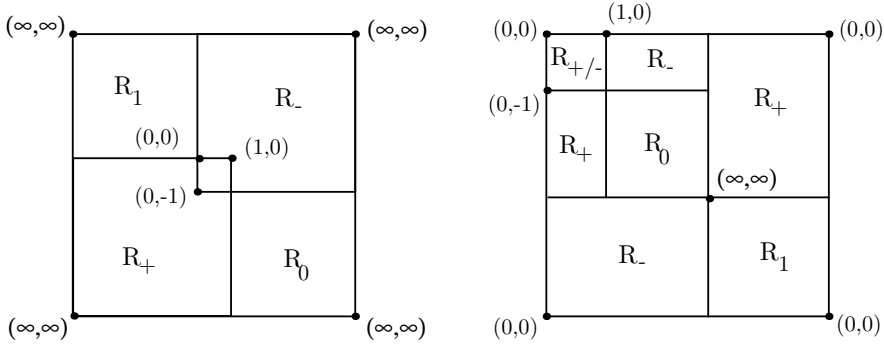}}
\caption{Partition of $\overline{\R^2}$ when $a<-1$.\label{partition}}
\end{figure}

Consider the covering of $\overline{\R^2}$ by 
closed rectangles:
\begin{eqnarray*}
R_+ & = & [-\infty,1] \times [-\infty,0] \\
R_- & = & [0,\infty] \times [-1,\infty] \\
R_0 & = & [1,\infty] \times [-\infty,-1] \\
R_1 & = & [-\infty,0]\times [0,\infty].
\end{eqnarray*}
Two views of this covering are pictured in Figure \ref{partition}; the right 
hand illustration is useful for visualizing the fixed point $(\infty,\infty)$, 
since in a small neighborhood of $(\infty,\infty)$, the action of $f$ is 
approximately $(s,t)\mapsto(t,s)$, which is a reflection about the diagonal 
$s=t$. 

\begin{prop} 
\label{thm1}
If $a<-1$ then the following hold:
\begin{itemize}
\item    
$f(R_+)\subset R_+$, and if $(x_0,y_0) \in R_+$ and $(x_1,y_1) = f(x_0,y_0)$, 
then
$$
\min\{x_1 - 1, y_1\} \leq \min\{x_0 - 1,y_0\} - 1.
$$
\item
$f^{-1}(R_-)\subset R_-$, and if $(x_0,y_0) \in R_-$ and $(x_{-1},y_{-1}) = 
f^{-1}(x_0,y_0)$, then
$$
\max\{x_{-1},y_{-1}+1\} \geq \max\{x_0,y_0+1\} + 1.
$$
\item $f(R_1)\cap \interior R_1\cap{\bf R}^2 = f^{-1}(R_1)\cap \interior R_1
\cap{\bf R}^2=\emptyset$.
\end{itemize} 
\end{prop}

\begin{figure}
\centerline{\epsfxsize3.8in\epsfbox{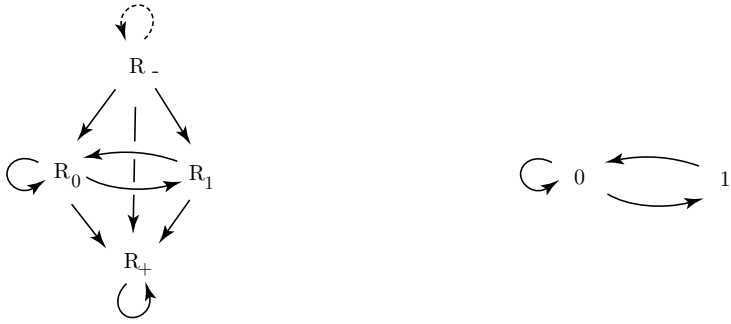}}
\caption{Graph of the Filtration \label{graph}}
\end{figure}
The proof of this proposition is elementary, and we leave its verification to the reader.  
Taken together, the four conclusions indicate that $f$ has the combinatorial 
behavior given by the graph on the left side of Figure \ref{graph}.  For 
instance, the arrow from $R_0$ to $R_1$ indicates
that $f(R_0)\cap R_1$ contains an open set.  The dashed arrow indicates 
that a point can remain within $\interior R_-$ for only finite positive time.  
The nontrivial recurrent part of this graph is given by the right hand 
side of Figure \ref{graph}.  If we reverse the arrows and move the dashed arrow to $R_+$,
then we obtain the graphs corresponding to $f^{-1}$.

Let us define the {\em stable set} $W^s(\infty,\infty)$ as 
the set of  points
$p$ such that $f^np\notin I(f)$ for all $n\ge0$, and
$\lim_{n\to+\infty}dist(f^np,(\infty,\infty))=0$. 

\begin{thm}  We have
$$W^s(\infty,\infty)=\bigcup_{n\ge0}f^{-n}R_+ -I_+.$$
In 
particular, if $f^np\in(R_0\cup R_1)\cap {\bf R}^2$ for all $n\in{\bf Z}$, then $f^np$ cannot 
approach $(\infty,\infty)$ in either forward or backward time. \end{thm}
\proof
If $p\in \bigcup_{n\ge0}f^{-n}R_+ -I_+$, then we have
$\lim_{n\to\infty}f^np=(\infty,\infty)$ by Proposition \ref{thm1}.  Conversely, suppose that
$p\notin I_+$, and 
$\lim_{n\to\infty} f^np=(\infty,\infty)$.    We will show that $p\in \bigcup_{n\ge0}f^{-n}R_+$. 
If not, then $f^np\in \interior(R_-\cup R_0\cup R_1)$ 
for all $n\ge0$; in particular, $f^np\in \R^2$.  By Proposition \ref{thm1} 
there can be only a finite interval $0\le j\le J$ for which $f^jp\in R_-$.  
Further, by Proposition \ref{thm1} once $f^{j_0}p\notin R_-$, we must have 
$f^jp\notin R_-$ for $j\ge j_0$. Thus we have $f^np\in R_0\cup R_1$ for 
$n\ge N$.

Without loss of generality, we may assume that $f^{2n}p\in R_0$ for all 
$n\ge0$.    We write $(x_n,y_n) = f^{2n}(p)$ and $m_n = y_n/x_n < 0$. 
Proposition \ref{approx} gives us
$$
m_{n+1} < m_n\left(1-\frac1{y_n}\right) 
$$
so that $|m_n|$ increases with $n$.  Therefore Proposition \ref{approx} also 
gives us that 
$$
y_{n+1} \geq y_n - C
$$
for some fixed $C$ and $n$ large enough.  From this, we deduce that 
$|y_n| \leq Cn $, that
$
\sum_{j=0}^\infty \frac{1}{y_n}
$
diverges, and that therefore
$$
\lim_{n\to\infty} m_n = m_0 \prod_{n=0}^\infty \frac{m_{n+1}}{m_n} = -\infty.
$$
So for $n$ large enough, we have $m_n = y_n/x_n < -1$ and 
$$
x_{n+1} = x_n + (a-1)x_n + m_n(a+1) + O(|x|^{-1} + |y|^{-1}) \leq x_n - 1 
$$
which contradicts the assumption that $x_n\to\infty$.  
\qed

\begin{cor}
If $p\in W^s(\infty,\infty)$, and if $f^jp\in R_0\cup R_1$ for all $j\ge0$, then we must have
$f^jp\notin{\bf R}^2$ for some $j\ge0$.  If in addition $p\in{\bf R}^2$, then we must have
$p\in f^{-n}\{x=1\}$ for some $n\ge0$.
\end{cor}

\section{Coding and Rectangles}
\label{coding}
\noindent Let us summarize some information about finite subshifts 
(see [KH] pages 176--181, and [LM]).  We use the 
symbol space $\cS=\{0,1\}^{\bf Z}$, which consists of bi-infinite sequences 
$w=\dots w_{-1}w_0\cdot w_1w_2\dots$, where $w_j\in\{0,1\}$ and the 
`$\cdot$' serves to locate the entry with subscript $0$.  Let
$\sigma:\cS\to\cS$ denote the shift operator given by $\sigma(w)=\ti w$, 
where $\ti w_j=w_{j+1}$.  We give $\cS$ the product space topology, which 
is generated by the finite cylinder sets 
$C(a_{-N}\dots a_N)\eqdef \{w\in\cS:w_j=a_j \text{\ for\ }-N\le j\le N\}$.
Thus $\cS$ is a compact space homeomorphic to a Cantor set.  Let us define
$\Sigma$ to be the subspace of $\cS$ consisting of all sequences (words) 
$w=(w_j)_{j\in{\bf Z}}$ such that the block ``11'' appears nowhere in 
$w$; alternatively, the symbol sequence $(w_j)$ may be generated by following the graph on the right hand 
side of Figure \ref{graph}.  We refer to $(\sigma,\Sigma)$ as the 
\emph{golden mean subshift}.

If $w\in\Sigma$ and if $j\le k$, we let $w[j,k]:=w_j\dots w_k$ denote the $[j,k]$
\emph{subword} of
$w$.  We refer to $[j,k]$ as the \emph{extent} of the word $w[j,k]$.  We let
$\Sigma^*$ denote all the subwords of elements $w\in\Sigma$.  If $w\in\Sigma^*$,
is a word of extent $[-n,m]$, with $-n\le0\le m$, we let $w^-=w[-n,0]$ denote the
$[-n,0]$ subword and $w^+=w[0,m]$  denote the $[0,m]$ subword. 

We say that a word $w$ is
\emph{admissible} if
$w\in\Sigma^*$ and if
$w$ has extent $[-n,m]$ with $0\le m,n\le\infty$.  We will only work with
admissible words in the rest of this paper, so we will use ``word'' to mean
``admissible word.'' 

 We let $\Sigma^+$ (resp. $\Sigma^-$) denote 
the sets of all $[0,\infty]$ words (resp. $[-\infty,0]$) words in $\Sigma^*$.   
We endow both spaces with the product topology.  The one-sided shift
$\sigma^+(w_0\cdot w_1\dots ) = w_1\cdot w_2\dots$ (resp.   
$\sigma^-(\dots w_{-1}w_0\cdot) = \dots w_{-2}w_{-1}\cdot$) gives a continuous 
self-map of $\Sigma^+$ (resp. $\Sigma^-$).

The incidence matrix corresponding to the graph on the right hand side of 
Figure \ref{graph} is
$
\left(\begin{matrix}
1&1 \\
1&0
\end{matrix}\right)
$, and the powers of this matrix satisfy
$
\left(\begin{matrix}
1&1 \\ 
1&0
\end{matrix}\right)^n=\left(\begin{matrix}F_{n+1}&F_n \\  
F_n& F_{n-1}\end{matrix}\right)
$, 
where $F_{-1}=1$, $F_0=0$, and $F_{n+1}=F_n+F_{n-1}$ denote the Fibonacci 
numbers.   For  $k,l\in\{0,1\}$, the $(k,l)$ entry of the $n$-th power of this 
matrix gives the number of words of length $n+1$ starting at $k$ 
and ending at $l$.  For example, if $-n\le0\le m$, then the number of $[-n,m]$
words that begin and end with `0' is $F_{n+m+1}$.

A sequence $w\in\Sigma$ satisfies $\sigma^nw=w$ if and only if $w_{n+j}=w_j$ 
for all $j\in{\bf Z}$.  Such a sequence is determined by a finite word 
$w_0\cdot w_1\dots w_n$ with $w_0=w_n$. We use the notation $\overline v:=\dots v\cdot vv\dots$ with $v=w_1\dots w_n$  to
express such periodic words. Since there two possibilities, 
$w_0=0$ or $w_0=1$, the number of such words is given by 
$$
\#\{w\in\Sigma:\sigma^nw=w\}=\text{ trace}\left(\begin{matrix}
1&1 \\ 
1&0
\end{matrix}\right)^n=F_{n+1}+F_{n-1}.$$

The topological entropy of $\sigma:\Sigma\to\Sigma$ is the logarithm of
the golden mean $\phi=(1+\sqrt5)/2$.   Every $\sigma$-invariant probability measure
on
$\Sigma$ has entropy less than or equal to $\log\phi$.  There is a unique
$\sigma$-invariant probability measure 
$\nu$ on $\Sigma$ with entropy equal to $\log\phi$.  This measure is 
given by averaging point masses over the periodic points 
$$
\nu=\lim_{n\to\infty}\frac 1{F_{n+1}+F_{n-1}} \sum_{\sigma^nw=w} \delta_w,
$$
where $\delta_w$ denotes the point mass at the point $w$.  The measure $\nu$ 
is also mixing.

There is a unique 
$\sigma^+$-invariant measure $\nu^+$ on $\Sigma^+$ with entropy $\log\phi$. 
The measures $\nu^\pm$ are \emph{balanced}:  if $E\subset\Omega^+$ is a 
measurable subset and $\sigma^+|_E$ is injective, then
$\sigma_+^*\nu^+|_E = \phi^{-1}\nu^+|_{\sigma^+(E)}$.  We may identify 
$\nu$ with $\nu^+\otimes\nu^-$ via the product structure
$$
\Sigma\ni v\mapsto\{(v^+,v^-)\in
\Sigma^+\times\Sigma^-:v^+[0]=v^-[0]\}.
$$

The stable manifold of a point $w\in\Sigma$ is
$$
W^s(w)=\{v\in\Sigma:v[n,\infty]=w[n,\infty]\text{ for some }n\in{\bf N}\}.
$$
We define the local stable manifold $$W^s_{loc}(w)=\{v\in\Sigma:v^+=w^+\},$$
which is the
cylinder $C(w^+)$ over the semi-infinite word $w^+$, and it follows that
$W^s(w)=\bigcup_{n\ge0}\sigma^{-n}W^s_{loc}(w)$.  It is evident that
$W^s_{loc}(\ti w)=W^s_{loc}(w)$ if and only if $\ti w^+=w^+$, so that the set 
of local stable manifolds is parametrized by the set $\Sigma^+$, and an 
individual local stable manifold $W^s_{loc}(w^+)$ is parametrized by
$\Sigma_-$.

 For a finite word 
$w$, we define the \emph{rectangle}
$$
R(w) \eqdef \overline{\bigcap_{k=-n}^mf^{-k}\,\interior R_{w_k}}.
$$
Since the interiors of $R_0$ and $R_1$ avoid $I(f^{-k})$ for all
$k\in\Z$, the definition of $f^{-k}$ is unambiguous.  When $w$ is a word of infinite
extent, we take $R(w)$ to be the intersection of all rectangles $R(w')$ 
corresponding to finite subwords $w'$ of $w$.   This definition of rectangle will
be shown in Theorem \ref{equivalentdef} to essentially coincide with two other
plausible definitions of rectangle. 

\begin{figure}
\centerline{\epsfxsize4.9in\epsfbox{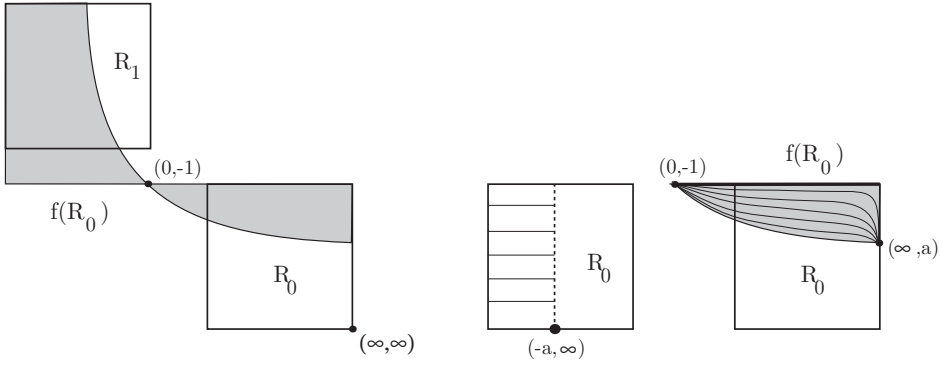}}
\caption{$R(00\cdot)$ \label{R00}}
\end{figure}

Let us consider an example when $a=-2$.  The left hand part of Figure 
\ref{R00} shows $R_0$, $R_1$, and $f(R_0)$.   The two squares on the right of 
Figure \ref{R00} show 
how $f$ maps part of $R_0$ onto $R({00\cdot})$, which is the intersection of the
shaded  region with $R_0$.  The dashed vertical segment inside $\{x=-a\}$ is
mapped  to the point $(0,-1)$.  The top segment in the boundary of $R_0$ is mapped
to  the curved portion of the boundary of $f(R_0)$.  The five horizontal segments 
$\{y=const, 1\le x\le -a\}$ in the interior of $R_0$, are mapped to the curves
connecting $(0,-1)$ and $(\infty,a)$ through the interior of $f(R_0)$.  The 
point of indeterminacy $(a,\infty)$ is mapped to the line $\{y=-1\}$, part of 
which forms the upper boundary of $R_0$.   This Figure will be discussed further in
connection with Figure 6.

The second item in the following Proposition shows why we assume that all words are
admissible.
\begin{prop} 
\label{interior}
Let $w$ be a finite $[-n,m]$ word, and let $R(w)$ be the corresponding rectangle.
Then
\begin{itemize}
\item $R(w) = \overline{\interior R(w)}$;
\item If $w$ is not admissible, then $R(w) =
\emptyset$;
\item $f^kR(w) = R(\sigma^k w)$ for all $-n\leq k \leq m$;  
\item $\interior R(w)\cap \cC(f^k) = \emptyset$ for $-n\leq k\leq m$.
\item $p \in \interior R(w)$ if and only if $f^kp\in \interior R_{w_k}$ for 
      $-n\leq k\leq m$.
\item For all $-n\leq k\leq m$, $f^k$ maps $\interior R(w)$ homeomorphically 
      onto $\interior R(\sigma^k w)$.
\end{itemize}
\end{prop}

\proof
The first conclusion follows because $R(w)$ is the closure of an open set.
The second follows from Proposition \ref{thm1}.  The third results from 
our convention for images of closed sets under $f$.  

To see that the fourth conclusion holds, suppose that it fails for some 
smallest $k>0$ (the case $k < 0$ is similar).  Then there exists 
$p\in\interior R(w)$ such that 
$f^{k-1}(p)\subset\cC(f)\cap \interior R_{w_{k-1}}$.  But this means that 
$f^{k-1}(p)\in\{x=-a\}$, and $f^k(p) = (0,-1)\notin R_{w_k}$, which 
conflicts with the assumption that 
$p\in R(w)\subset {f^{-k}\,R_{w_k}}$.  Hence the
fourth conclusion is true.

We also have that $\interior R(w)\subset \interior R_{w_0}$ avoids the 
indeterminacy set of every iterate of $f$.  So in light of the fourth
conclusion,the restriction $f^k|_{\interior R(w)}$ is a homeomorphism
for each $-n\leq k\leq m$.  The fifth and sixth conclusions follow
immediately.
\qed

We will call a connected arc $\gamma\subset R_0$ an \emph{s-arc} if it joins 
the boundary components $\{y=-1\}$ and $\{y=-\infty\}$ and a \emph{u-arc} if 
it joins $\{x=1\}$ and $\{x=\infty\}$.  Similarly, we call $\gamma\subset R_1$
an s-arc if it joins $\{x=0\}$ and $\{x=\infty\}$ and a u-arc if it joins 
$\{y=0\}$ and $\{y=\infty\}$.  In any case, let us call $\gamma$ \emph{proper}
if only its endpoints lie on the boundaries of $R_0$ and $R_1$ and these 
points are not corners.

\begin{prop}
\label{allowables}
Let $\gamma$ be a proper arc in $R_0\cup R_1$.  
\begin{itemize}
\item If $\gamma\subset R_0$ is a u-arc, then $f(\gamma)$ contains 
      proper u-arcs in $ R_0$ and $R_1$.
\item If $\gamma\subset R_1$ is a u-arc, then $f(\gamma)$ contains a
      proper u-arc in $R_0$.
\item If $\gamma\subset R_0$ is an s-arc, then $f^{-1}(\gamma)$ contains
      proper s-arcs in $R_0$ and $R_1$.      
\item If $\gamma\subset R_1$ is an s-arc, then $f^{-1}(\gamma)$ contains a
      proper s-arc in $R_0$.
\end{itemize}
\end{prop}

\proof
Suppose that $\gamma\subset R_0$ is a proper u-arc.
By the intermediate value theorem, there is a smallest subarc 
$\alpha\subset\gamma$ beginning with the left endpoint of $\gamma$ and ending 
on $\{x=-a\}$.  The map $f$ sends the left endpoint of $\alpha$ to the line 
$(a,\infty)$ and the right endpoint of $\alpha$ to $(0,-1)$.  One sees easily 
that the intervening points are mapped into the region 
$(-\infty,0)\times(-1,\infty)$.  Hence $f(\alpha)\cap R_1$ contains a 
proper u-arc.  Likewise, there is a smallest subarc 
$\beta\subset\gamma$ beginning on $x=-a$ and ending with the left endpoint 
of $\gamma$ on $x=\infty$.  The left and right endpoints of $\gamma$ are sent 
to $(0,-1)$ and $(b,\infty)$, respectively, where $b<-1$.  The intervening 
points are all mapped below $y=-1$, so that $f(\beta)\cap R_0$ must contain a
proper u-arc.  The first assertion is now proved.  The proof the second assertion
is similar.  The last two assertions follow from the reversibility of $f$. 
\qed

\begin{figure}
\centerline{\epsfxsize1.8in\epsfbox{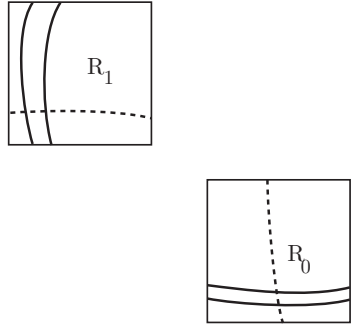}}
\caption {Proper unstable (solid) and stable (dashed) arcs.\label{arcs}}
\end{figure}

\begin{prop}
\label{lowerbd}
Let $w = w[-n,0]$ be a finite word, and let $L$ be a horizontal or vertical
line that intersects $\interior R_{w_{-n}}$ in a proper u-arc.  Then $f^n(L)\cap
\interior R(w)$  contains a proper u-arc $\gamma$.  Likewise, if $w = w[0,m]$ is a
finite  admissible word and $L$ meets $R_{w_m}$ in a proper s-arc.  Then 
$f^{-m}(L)\cap R(w)$ contains a proper s-arc $\gamma$.  
\end{prop}

\proof
We work by induction, considering only the case $w = w[-n,0]$.
By hypothesis $L\cap R_{w_{-n}}$ is a proper u-arc $\gamma_{-n}$.  Suppose 
that for $-n \leq -j < -k \leq 0$ we have proper u-arcs 
$\gamma_j\subset R_{w_{-j}}$ satisfying $\gamma_{j-1}\subset f \gamma_j$.
Then since $w$ is admissible, Proposition \ref{allowables} gives us a proper 
u-arc $\gamma_k\subset f \gamma_{k+1}\cap R_{w_{k+1}}$.  Hence $\gamma_j$ 
exists for $0\leq j\leq n$.  Moreover, if $p\in\gamma_0$ is not an endpoint, 
then neither is $f^{-j}(p)\in\gamma_j$ for any $j$.  Hence the portion 
$f^{-n}(p),\dots , p$ of the orbit of $p$ lies entirely in 
$\interior R_0\cup \interior R_1$.  It follows that $p\in \interior R(w)$.  
Hence $\gamma = \gamma_0$ is the arc we are seeking.
\qed

\begin{thm}
\label{nonempty}
 $R(w)$ is nonempty.
\end{thm}

\proof
Since rectangles corresponding to finite words are compact, and rectangles 
corresponding to infinite words are decreasing intersections of these, it is 
enough to prove the proposition for a finite word $w$ of arbitrary extent 
$[-n,m]$.  
To do this, we apply the previous proposition to obtain a proper u-arc 
$\gamma^-\subset R(w^-)$ and a proper s-arc $\gamma^+\subset R(w^+)$.  Then 
$\gamma^-,\gamma^+\subset R_{w_0}$ must meet at some point $p$ that is not 
an endpoint of either arc.  It follows that  
$
p\in \interior R(w^+)\cap \interior R(w^-) = \interior R(w),
$
where the equality holds by the definition of $R(w)$ and the fifth 
conclusion of Proposition
\ref{interior}.
\qed

\section{Complexification and Intersection Theory}
\label{intersection}
\noindent  We will work with the complexification $\tilde f$ of $f$; for
simplicity, we drop the tilde.  The homology
group
$H_2(\pp,\C)$ is two dimensional.  Let us fix generators: $\gamma_1$, which
corresponds to the horizontal (complex) line $\cp^1\times\{y_0\}$, and
$\gamma_2$, which corresponds to the vertical (complex) line
$\{x_0\}\times\cp^1$.  Let
$\{\gamma_1^*,\gamma_2^*\}$ denote the basis of $H^2(\pp,\C)$ which is dual
to
$\{\gamma_1,\gamma_2\}$.  We may represent $\gamma_1^*$ as the (1,1) form
$\frac{i}{2\pi}(1+|x|^2)^{-2}dx\wedge d\bar x$ and $\gamma_2^*$ as
$\frac{i}{2\pi}(1+|y|^2)^{-2}dy\wedge d\bar y$.  The homology classes of the 
preimages under $f$ are:
$f^{-1}\gamma_2\sim \gamma_1+\gamma_2$ and $f^{-1}\gamma_1=\gamma_2$.  Thus
the induced pullback map $f^*$ on the cohomology group $H^2(\pp)$ is
given with respect to the basis $\gamma^*_1$, $\gamma_2^*$ as
$$
\label{matrix}
f^*=\left(\begin{matrix}
0&1 \\ 
1&1
\end{matrix}\right).
$$
In particular, the largest eigenvalue of $f^*$ is the golden mean $\phi$.  
By [DF], the pullback of $f^n$ on cohomology, $(f^n)^*$, coincides
with $(f^*)^n$ if $a$ is not exceptional.   Then the powers of $f^*$ are
given by the Fibonacci numbers:
$$
\left(\begin{matrix}
0&1 \\
1&1
\end{matrix}\right)^n
=
\left(\begin{matrix}
F_{n-1} & F_n \\ 
F_n & F_{n+1}
\end{matrix}\right).
$$
So we have
$f^{n*}\gamma_1^*=F_{n-1}\gamma^*_1+F_n\gamma_2^*$ and $f^{n*}\gamma^*_2 =
F_{n}\gamma_1^*+ F_{n+1}\gamma_2^*$. 

Complex algebraic curves $V$ and $W$ in ${\bf P}^1\times{\bf P}^1$ define
cohomology classes $\{V\}=n_1\gamma^*_1+n_2\gamma^*_2$ and
$\{W\}=m_1\gamma^*_1+m_2\gamma^*_2$ in $H^2$.  We use the notation
$V\sim[n_1,n_2]$ and $W\sim[m_1,m_2]$. The intersection product on $H^2$ is
defined as
$$V\cdot W=n_1m_2 + n_2m_1.$$  Thus the
\emph{intersection form} 
$$\cdot: H^{2}(\pp)\times H^{2}(\pp)\to \C$$
is a quadratic form whose matrix is
$$
\left(\begin{matrix}
0&1 \\
1&0
\end{matrix}\right) 
$$
with respect to the basis $\gamma_1^*,\gamma_2^*$.  A basic result of
intersection theory (see [Fu]) is that if all points of $V\cap W$ are
isolated, then the intersection product $V\cdot W$ is equal to the number of
intersection points
$V\cap W$ counted with multiplicity.  Since the curves are complex, the
multiplicity of each isolated intersection point is an integer $\ge1$.

The pushforward $f_* = (f^{-1})^*$ of $f$ acting on $H^{2}(\pp)$ is just the
adjoint
$$
\left(\begin{matrix}
1&1 \\
1&0
\end{matrix}\right) 
$$
of $f^*$ with respect to this form.  Thus
$f_*^m\gamma_1^*=F_{m+1}\gamma_1^*+F_m\gamma_2^*$, and $f_*^m\gamma_2^* =
F_m\gamma_1^*+F_{m-1}\gamma_2^*$.  In the sequel we will use this to compute
the number of intersections between $f^{n*}\gamma_i$ and $f^m_*\gamma_j$. 
For instance:
$$
f^{n*}\gamma_2\cdot f^m_*\gamma_1 = F_{n+1}F_{m+1} + F_nF_m.
$$

\begin{lem} For each $n\in\N$, the periodic points of period $n$ for $f$ 
are isolated.
\end{lem}

\proof
The alternative is that $f^n$ fixes some algebraic curve $V\subset\pp$ pointwise.  
But $V$ is homologous to a positive linear combination of horizontal and vertical lines
and must therefore intersect the line $\{y=x-1\}$ somewhere.  Since $f$ acts by
translation on this line with $(\infty,\infty)$ as the sole fixed point, we see that
$V$ must contain the point $(\infty,\infty)$.  This contradicts Proposition \ref{isolated}
below.
\qed

Given the Lemma, let us describe how to count the
periodic points of the complexification 
of $f$.  Let $\Gamma_{f^n}$ denote the graph of $f^n$ as a subvariety of
$(\pp)\times(\pp)$.  The periodic points of $f$ are given by intersecting
$\Gamma_{f^n}$ with the diagonal $\Delta\subset(\pp)\times(\pp)$. Thus we
have:
$$
\{p\in\pp:f^np=p\} =\Delta\cap \Gamma_{f^n}.
$$
Since $f^{-n}(I(f))\cap I(f) = \emptyset$ for every $n\in\N$,
it follows from the Lefschetz Fixed Point Formula [Fu]
that
$$
\#\{p\in\pp:f^np=p\}=\text{trace}(f^*),
$$
where the periodic points are counted according to their multiplicities, and
the trace refers to the action of $f^*$ on cohomology in all dimensions. 
The total cohomology is given by $H^*(\pp)=H^0+H^2+H^4$.  Now
$f^*$ acts as the identity on 
$H^0$ and $H^4$, both of which have dimension 1, so we may evaluate the total
trace to obtain
$$
\#\{p\in\pp:f^np=p\}=F_{n+1} + F_{n-1}+2.
$$

We have seen that $(\infty,\infty)$ is the only fixed point of $\pp-{\bf C}^2$.
Thus we have
$$
\#\{p\in{\bf C}^2:f^np=p\}=F_{n+1} + F_{n-1}+2-m_{(\infty,\infty)},
$$
where $m_{(\infty,\infty)}$ denotes the multiplicity of $(\infty,\infty)$ as
a fixed point of $f^n$, which is defined as the multiplicity of the
intersection of $\Delta$ and $\Gamma_{f^n}$ at $(\infty,\infty)$.

\begin{prop}
\label{isolated}
For every $n\in{\bf Z}$, $n\ne 0$, the point $(\infty,\infty)$ is an isolated
fixed point of
$f^n$ with multiplicity at least two.  When $n$ is even,
$m_{(\infty,\infty)}=4$.
\end{prop}

\proof
Writing $f$ with respect to the coordinates
$(\xi,\eta) = (1/x,1/y)$ and employing Proposition \ref{approx} gives
$$
f^2(\xi,\eta) = (\xi,\eta) + Q(\xi,\eta) + O(\norm{(\xi,\eta)}^3),
$$
where $Q(\xi,\eta) = (\xi^2(1-a) - \xi\eta(1+a),\eta^2(1-a) - \eta\xi(1+a))$ 
is a non-degenerate homogeneous map of degree 2.  
Therefore , 
$$ 
f^{2n}(\xi,\eta) = (\xi,\eta) + nQ(\xi,\eta) + O(\norm{(\xi,\eta)}^3)
$$
for every $n\in\Z$.  In particular, there exists $C>0$ and 
$\epsilon = \epsilon(n)$
$$
\norm{f^{2n}(\xi,\eta) - (\xi,\eta)} \geq C\norm{(\xi,\eta)}^2
$$
for every $\norm{(\xi,\eta)} < \epsilon$.  So $(\xi,\eta)=(0,0)$ is isolated
as a fixed point of $f^{2n}$.  Fixed points of $f^n$ are also
fixed by $f^{2n}$, so $(\xi,\eta)=(0,0)$ is isolated as a fixed point of
any iterate of $f$.

The multiplicity of an isolated fixed point of $f^n$ is greater than one 
exactly when one of the eigenvectors of $Df^n$ has eigenvalue one.
One can easily check that $(1,1)$ is such an eigenvector for 
$Df^n_{(\xi,\eta)=(0,0)}$.  To compute the exact multiplicity in the even 
case, we need to compute the multiplicity of $(0,0)$ as a solution of 
$f^{2n}(\xi,\eta) - (\xi,\eta) = 
nQ(\xi,\eta) + O(\norm{(\xi,\eta}^3) = (0,0)$.  Because
$Q$ is non-degenerate and quadratic, it follows that the multiplicity
is four.
\qed

It will be seen in \S6 that $m_{(\infty,\infty)}=2$ when $n$ is odd.

\section{Structure of rectangles}
\label{productstructure}
\noindent   Here we study the structure of rectangles $R(w)$ for finite $w$. 
Essentially, we continue \S2, now incorporating complex intersection theory.  The
first result (Theorem \ref{interior2}) is that if $w$ is finite, then the interior
of $R(w)$ has a canonical product structure.  This product structure extends to
points of $R(w)\cap{\bf R}^2$ but degenerates at $R(w)-{\bf R}^2$.  Problems with
points at infinity lead us to consider the special case of words $w$
which are ``alternating.'' The possibilities for (nonempty)
$R(w)\cap{\bf R}^2$ are given in Theorem
\ref{whatsinr} and Corollary \ref{infpts}.  Theorem \ref{whatsinr} then leads to
characterizations of $R(w)$ (Theorem 4.9) and of
$\Omega:=\bigcup_{w\in\Sigma}R(w)$ (Theorem 4.11).

\begin{thm} 
\label{itsanarc}
Let $w_0,w_{-n}\in\{0,1\}$ be given, and let $W$ denote the set of 
$[-n,0]$ words $w$ beginning with $w_{-n}$ and ending with $w_0$.  Let $L$ be 
a horizontal or vertical complex line that meets $R_{w_{-n}}$ in a proper 
u-arc.  Then
$$
f^n L\cap R_{w_0} = \bigcup_{w\in W} R(w)\cap f^n L,
$$ 
and for each $w\in W$, $R(w)\cap f^n L$ is a proper u-arc.  If $w_0 = 0$
(respectively $w_0 = 1$) then $R(w)\cap f^n L$ can be expressed as the
graph of a function over the $x$-axis (respectively, $y$-axis). 
\end{thm}

\proof
We treat only the representative case $w_{-n} = w_0 = 0$ 
(in particular $L$ is horizontal); the 
other cases are similar.  By Proposition \ref{lowerbd} $f^n L\cap\interior R(w)$
contains a proper  u-arc $\gamma_w$ for each $w\in W$.  By Proposition 2.1, we know
that
$\gamma_w\cap\gamma_{\ti w} = \emptyset$ for distinct words 
$w,\ti w\in W$.  Moreover, the definition of a u-arc in $R_0$ implies that 
$\gamma_w$ intersects the vertical line $\{x=x_0\}$, for each 
$x_0\in [1,\infty]$ and each $w\in W$.  Hence if $x_0\in (1,\infty)$, there 
are at least $F_{n+1}$ distinct intersections between $\{x=x_0\}$ and $f^nL$.  On
the  other hand, as complex curves $\{x=x_0\}\cdot f^nL = F_{n+1}$, too.  Hence 
there are no further intersections between $\{x=x_0\}$ and $f^nL$.  Since 
vertical lines foliate $R_0$, we see that
$f^n L\cap R_0 = \bigcup_{w\in W}\gamma_w$ and that $f^n L\cap R(w)$ is 
a graph over the $x$-axis, as desired.
\qed

Let us define projections $\pi^{s/u}_j:\pp\to\cp^1$ for $j=0,1$ according to 
the formulas
$$
\pi^s_0(x,y) = \pi^u_1(x,y) = x,\quad \pi^u_0(x,y) = \pi^s_1(x,y) = y.
$$
We choose intervals 
$$
T^u_0 = [-\infty,-1], \quad
T^s_0 = [1,\infty], \quad
T^s_1 = [0,\infty], \quad
T^u_1 = [-\infty,0]
$$ 
so that for $j=0,1$, the map $\pi_j \eqdef (\pi^s_j,\pi^u_j):\pp\to\pp$ is a 
biholomorphism that restricts to a homeomorphism from $R_j$ onto 
$T^s_j\times T^u_j$.  

More generally, if $w = w[-n,m]$ is a finite word, let us define 
$\pi^{s/u}_w:\pp\to\cp^1$ by 
$$
\pi^s_w = \pi^s_{w_m}\circ f^m, \quad \pi^u_w = \pi^u_{w_{-n}}\circ f^{-n},
$$
and set $T^s_w = T^s_{w_m}$, $T^u_w = T^u_{w_{-n}}$. This gives us a
meromorphic map
$\pi_w := (\pi^s_w,\pi^u_w):\pp\to\pp$ whose restriction 
maps $\interior R(w)$ to $\interior T^s_w\times T^u_w$.   Note that
$\pi_w$ depends only on the first and last digits of $w$.  Clearly,
$$
\pi_{\sigma^k w}\circ f^k = \pi_w
$$
for $-n\leq k \leq m$.  Furthermore, for each $t\in T^u_{w_{-n}}$ we have from
Theorem \ref{itsanarc} that $R(w)\cap (\pi^u_w)^{-1}t$ is the intersection
of
$R(w)$  with a proper u-arc.  Likewise for each $t\in T^s_{w_m}$, the
fiber 
$R(w)\cap (\pi^s_w)^{-1}t$ is the intersection of $R(w)$ with a proper
s-arc.  We refer  to these fibers as \emph{canonical  s/u-arcs of
$R(w)$.}  Two properties of canonical arcs follow immediately from their
definition.  We  state them only for u-arcs.
\begin{itemize}
\item If $\gamma$ is a canonical u-arc of $R(w)$, and if $\sigma w$ is
admissible, then  
      $f\gamma$ is a canonical u-arc of $R(\sigma w)$.
\item If $\ti w$ extends $w$ to the right, then 
      the canonical u-arcs of $R(\ti w)$ are sub-arcs of the canonical
u-arcs of $R(w)$.
\end{itemize}
Where context makes things clear, we will drop the subscripts from 
$\pi^u_w,T^u_{w_m}$, etc.

We say that a curve $V$ belongs to the \emph{exceptional
locus}  of $\pi_w$ if $\pi_w(V)$ is a point.

\begin{prop}
\label{exceptional}
For any finite word $w = w[-n,m]$, the indeterminacy locus of 
$\pi_w$ is contained in $I(f^{-n})\cup I(f^m)$.  The intersection between 
$R(w)$ and the exceptional locus of $\pi_w$ is contained in 
$\pi_w^{-1}(\infty,\infty)\cap (\{x=\infty\}\cup\{y=\infty\})$.
\end{prop}

\proof
The assertion about the indeterminacy set is clear from the definition of 
$\pi$.  To see that the claim about the exceptional locus of $\pi$ is true, 
note that $\pi(R(w)-(I(f^{-n})\cup I(f^m))) \subset T^s\times T^u$.  So fix
$(x_0,y_0)\in T^s\times T^u$.  Then any overlap between $\pi^{-1}(x_0,y_0)$ 
and the exceptional set of $\pi$ is a common component of the complex curves 
$(\pi^s)^{-1}x_0$ and $(\pi^u)^{-1}y_0$.

Since $(\pi^s_{w_m})^{-1}x_0$ is a line in $R_0\cup R_1\cup R_-$, it follows 
from Proposition \ref{thm1} that any irreducible component of
$(\pi^s)^{-1}x_0 = f^{-m}(\pi^s_{w_m})^{-1}x_0$ not equal to $\{x=\infty\}$ or 
$\{y=\infty\}$ must contain points in $R_-- R_+$.  Similarly, any 
non-infinite irreducible component of $(\pi^u)^{-1}y_0$ contains points in 
$R_+- R_-$.  Therefore the only candidates for a common irreducible 
component $V$ of $(\pi^s)^{-1}x_0$ and $(\pi^u)^{-1}y_0$ are $\{x=\infty\}$ 
and $\{y=\infty\}$, and it follows that $f^k(V)$ is  $\{x=\infty\}$ or 
$\{y=\infty\}$ for all $k\in\Z$.  Since $\pi(V) = (\pi^u(V),\pi^s(V))$ is a 
single point, we must have $\pi(V) = (\infty,\infty)$.
\qed

As the next theorem shows, the restriction of $\pi_w$ to $\interior R(w)$
defines a product structure.

\begin{thm}
\label{interior2}
If $w$ is a finite word, then $\pi_w$ maps $\interior R(w)$ 
homeomorphically onto $\interior T^s\times \interior T^u$.  More generally, 
$\pi$ is injective on $R(w)\cap\R^2$.  So if $R(w)\subset\R^2$, 
then $\pi_w:R(w)\to T^s\times T^u$ is a homeomorphism.  
\end{thm}

\proof 
Let $[-n,m]$ be the extent of $w$.  Since 
$f^{-n}\,\interior R(w)\subset \interior R_{w_{-n}}$ and 
$f^m\,\interior R(w)\subset R_{w_m}$, we have that $\pi$ maps $\interior R(w)$ 
into $\interior (T^s\times T^u)$.

Fix a point $(x_0,y_0)\in \interior(T^s\times T^u)$.  Then the canonical u-arc 
$(\pi^u)^{-1}y_0\cap R(w^-)$ must meet the canonical s-arc 
$(\pi^s)^{-1}x_0\cap R(w^+)$ at some point 
$p_w\in\interior R(w^+)\cap \interior R(w^-) = \interior R(w)$.  So 
$\pi(p_w) = (x_0,y_0)$.  That is, $\pi$ is surjective.

On the other hand, $\pi$ depends only on the extent $[-n,m]$ of $w$ and the 
first and last digits $w_{-n},w_m$.  Hence our argument produces a distinct 
preimage of $(x_0,y_0)$ in $\interior R(\ti w)$ for every $[-n,m]$ word 
$\ti w$ whose first and last digits agree with those of $w$.  It is 
straightforward to verify that regardless of $w_{-n}$ and $w_m$, the number of 
such words is exactly the same as the intersection number of $(\pi^u)^{-1}y_0$ 
and $(\pi^s)^{-1}x_0$ treated as complex curves.  Therefore (as we argued in 
Theorem \ref{itsanarc}) there is exactly one preimage of $(x_0,y_0)$ in 
$\interior R(\ti w)$ for each $\ti w$ and no other preimages in $\pp$.  
In particular, $\pi$ is injective on $R(w)$.

We also obtain that fibers of the meromorphic map $\pi$ are discrete over 
points in $\interior(T^s\times T^u)$.  They are therefore discrete over
a neighborhood $U$ of $\interior(T^s\times T^u)$ in $\pp$.  Taking $U$ 
small enough, we see that $\pi^{-1}U$ is a disjoint union of connected 
components $U(\ti w)\subset R(\ti w)$ for each $[-n,m]$ word $\ti w$ whose 
first and last digits agree with those of $w$.  Because the number of such 
words is exactly the topological degree of $\pi$, we see that $\pi$ sends 
$U(\ti w)$ holomorphically and injectively onto 
$U$ for each $\ti w$.  We conclude that $\pi$ restricts to a homeomorphism 
from $\interior R(w)$ onto $\interior(T^s\times T^u)$.

By continuity, we must have that fibers $\pi^{-1}(x_0,y_0)\cap R(w)$ of the 
restricted map are connected even when $(x_0,y_0)\in \partial(T^s\times T^u)$. 
In addition, a point $p\in R(w)\cap {\bf R}^2\cap \pi^{-1}(x_0,y_0)$ 
in the finite part of a fiber must be isolated by Proposition
\ref{exceptional}.  Hence
$\pi$ is  injective on $R(w)\cap{\bf R}^2$.  In particular if $R(w)\subset \R^2$,
then
$\pi$ maps 
$R(w)$ homeomorphically onto its image in $T^s\times T^u$.  The image is 
compact and contains $\interior( T^s\times  T^u)$ as a dense subset,
so it is in fact  equal to $T^s\times T^u$.
\qed

We single out an observation from the proof of the preceding theorem as a separate
result.
\begin{prop}
Given $a,b\in\{0,1\}$ let $\pi:=(\pi^s_af^m,\pi_bf^{-n})$.  Then
$$\pi^{-1}(\interior (T^s_a\times T^u_b)) =\bigcup_w \interior R(w),$$ where the
union is taken over all $[-n,m]$ words $w$ with $w_{-n}=b$ and $w_m=a$.
\end{prop}

\begin{cor}
\label{connected}
$R(w)$ is connected.
\end{cor}

\proof
If $w$ is finite, then $R(w) = \overline{\interior R(w)}$ is connected by 
Theorem 4.1.  If $w$ is infinite, it is a decreasing intersection of 
compact, connected sets and must also be connected.
\qed

There are exactly two words in $\Sigma$, $\overline{01}$ and $\overline{10}$, in
which which the digits alternate.  We call a word {\em alternating} if it is a
subword of one of these.   We now describe the connection between rectangles that
contain points of indeterminacy and (partially) alternating words.  First an
elementary observation.

\begin{lem}
Suppose $j,k\in\{0,1\}$ are not both $1$, that $p\in R_j-\R^2 - I(f)$, that 
$f(p)\in R_k - \R^2$, and that neither point is a corner of $R_0$ or $R_1$.  
Then for any small neighborhood $U\ni p$, we have $f(U)$ is a neighborhood of 
$p$ and that $f(U\cap R_j) = f(U)\cap R_k$.
\end{lem}

\proof Since $p\in R_j-{\bf R}^2$ and $p\ne(-a,\infty)$, it follows that
$f$ is a local diffeomorphism at $p$.  The Lemma follows because $({\bf
R}\times\{\infty\})\cup (\{\infty\}\times{\bf R})$ is an invariant set,
and locally near $p$, $f$ maps $R_j$ to $R_k$. \qed

\begin{thm}
\label{whatsinr}
Let $w$ be a $[-n,m]$ word.  
\begin{itemize} 
\item $(\infty,\infty)$ belongs to $R(w)$ if and only if $w$ is alternating. 
\item If $p\notin I(f^{-n})\cup I(f^m)$ is not $(\infty,\infty)$, then 
      $p\in R(w)$ if and only if $f^k(p)\in {R_{w_k}}$ for 
      $-n\leq k\leq m$.  
\item $p\in I(f^m)$ belongs to $R(w)$ if and only if 
      $f^kp=(-a,\infty)$ for some $0\leq k< m$, and $w[-n,k]$ is 
      alternating, but $w[-n,k+1]$ is not.  
\item $p\in I(f^{-n})$ belongs to $R(w)$ if and only if 
      $f^{-k}p=(\infty,a)$ for some $0\leq k < n$, and $w[-k,m]$ is
alternating, 
      but $w[-k-1,m]$ is not.
\end{itemize}
If $w$ and $\ti w$ are distinct $[-n,m]$ words, then 
$R(w)\cap R(\ti w)$ contains at most one point, which is in $I_+\cup I_-$, and
consequently
$R(w)\cap R(\tilde w)\cap {\bf R}^2=\emptyset$.
\end{thm}

\proof
We need only establish each conclusion in the case where $w$ is finite. 
The first conclusion holds because points $q\in R_0$ near $(\infty,\infty)$
map to points 
$f(q)\in R_1$ near $(\infty,\infty)$. 

For the second assertion, let us suppose first that $q\in\interior R(w)$.  Then by
the fifth item in Proposition 2.1, 
$f^jq\in R_{w_j}$ for all $-n\le j\le m$.  Thus
$\Rightarrow$ holds in this case.  Now $f^j$ is continuous at $p$ since
$p\notin I(f^{-n})\cup I(f^m)\cup\{(\infty,\infty)\}$.  If we approximate
$p$ by $q\in\bigcap_{j=-n}^m f^{-j}\interior R_{w_j}$, we see that
$\Rightarrow$ holds by continuity.

For the case $\Leftarrow$ in the second assertion, we consider first the
case $p\in{\bf R}^2$.  By Theorem 4.2, $\pi_w$ is holomorphic and open at
$p$, and we have observed that $\pi_w(p)\in T^s_{w_m}\times
T^u_{w_{-n}}$.  Let us choose points $q_j\in\interior( T^s_{w_m}\times
T^u_{w_{-n}})$ converging to $\pi_w(p)$, and let us choose preimages
$p_j\in\pi^{-1}_w(q_j)$ which converge to $p$ as $j\to\infty$.  By Theorem
4.4, we may pass to a subsequence so that $p_j\in\interior R(\tilde w)$ for
some word $\tilde w$ starting with $w_{-n}$ and ending with $w_m$.  Since
$R(\tilde w)$ is closed, we must have $p\in R(\tilde w)$.  It follows by
the $\Rightarrow$ case that $f^jp\in R_{\tilde w_j}$.  Thus $\tilde w=w$.

If $p\notin{\bf R}^2$, then $f^jp\notin{\bf R}^2$.  The $\Leftarrow$ part
of the second assertion then follows by the preceding Lemma.

Next we prove the third assertion.  Since $p\in I(f^m)$, there exists a
unique $k$, $0\le k<m$ such that $f^kp\in I(f)$.  In particular $f^kp\in
I(f)\cap (R_0\cup R_1)=(-a,\infty)$.  Now if $p\in R(w)$, we may choose a
sequence $p_l\to p$ such that $f^jp_l\in int(R_{w_j})$ for all $-n\le j\le
m$.  By continuity, $f^j_{p_l}\to f^jp$ for all $-n\le j\le k$.  Note that
$f^kp\in R_0$, $f^{k-1}p\in R_1$, etc., so by the second item in this
Proposition, $p\in R(w[-n,k])$, and $w[-n,k]$ is alternating.  Although
$f$ is not continuous at $(-a,\infty)$, the points $f^{k+1}p_l$ must
accumulate only on $\dot f(-a,\infty)\cap (R_0\cup R_1)\subset
R_0\cap\{y=-1\}$.  Thus $f^{k+1}p_l\in R_0$, and so $w_k=w_{k+1}=0$.  Thus
$w[-n,k+1]$ is not alternating.  This proves the $\Rightarrow$ portion of
the third assertion.

To prove the $\Leftarrow$ part of the third assertion, define $v=0\cdot
v_1\dots v_{n-k-1}$ by setting $v_\ell = w_{k+1+\ell}$.  By Proposition
\ref{lowerbd}, $f^{m-k-1}\{y=-1\}\cap R(v)$ contains a proper u-arc $\gamma$. 
Choose
$q\in f^{-m+k+1}(\gamma\cap{\bf R}^2)$.  Since $q\in\{y=-1\}\cap{\bf
R}^2$, we have $f^{-1}q=(-a,\infty)$, and so $f^{-j}q\notin I(f^{-1})$ for
all $j\ge0$.  Thus we may apply the second item of this Proposition to
conclude that $q\in R(\tilde v)$ with $\tilde v=\overline{10} v$; we see
that we may concatenate the $\overline{10}$ on the left of $v$ because
$f^{-j}q\in R_{w_{-j}}$ for alternating symbols $w_{-j}$.  By Proposition
2.1, we have $p\in f^{-k-1}q\in R(\sigma^{-k-1}\tilde v)\subset R(w)$.

The proof of the fourth assertion is similar.  \qed

We can now specialize Theorem \ref{whatsinr} to the case of infinite points.

\begin{cor}
\label{infpts}
For any word $w$, we have the following possibilities:  
\begin{itemize}

\item If $w$ is alternating, then $R(w)\cap\{x=\infty\}$ and 
      $R(w)\cap\{y=\infty\}$ are (possibly degenerate) intervals containing
$(\infty,\infty)$.
\item If $w$ is not alternating, but there exists $n \leq k < m$ such that 
      $w[-n,k]$ and $w[k+1,m]$ are alternating, then $R(w) -\R^2$ is the 
      interval $f^k(E)$, where
$$
      E\eqdef \{(x,\infty):1 \le x \le -a\}.
$$
This case corresponds to $w$ being a subword of a
translate of $\overline{10}\cdot\overline{01}$. 
\item If neither $w^+$ nor $w^-$ is alternating, then $R(w)\subset\R^2$.
\item Otherwise, $R(w) - \R^2$ contains exactly one point, and this 
      point belongs to $I(f^{-n})\cup I(f^m)$.
\end{itemize} 
\end{cor}

\proof  By Theorem \ref{whatsinr}, we know that $p$ belongs to a rectangle
$R(w)$ if and only if it is part of an orbit with $w$ as its itinerary. 
The first assertion of this Corollary is immediate.  The third and fourth
items of Theorem \ref{whatsinr} assert that $R(w)$ can contain at most one element
of
$I(f^{-n})\cup I(f^m)$.  Further, if there is such an element, then the
Theorem says that $w^+$ or $w^-$ is alternating.

It remains to consider points $p\in R(w)-({\bf R}^2\cup I(f^{-n})\cup
I(f^m))$.  The orbit of such a point will alternate between $R_0$ and
$R_1$ unless there is a $j$ with $f^jp\in E$.  In this case we have
$f^jp,f^{j+1}p\in R_0$.  Note that this can happen for at most one $j$. 
This completes the proof.  \qed

\begin{figure}
\centerline{\epsfxsize3.2in\epsfbox{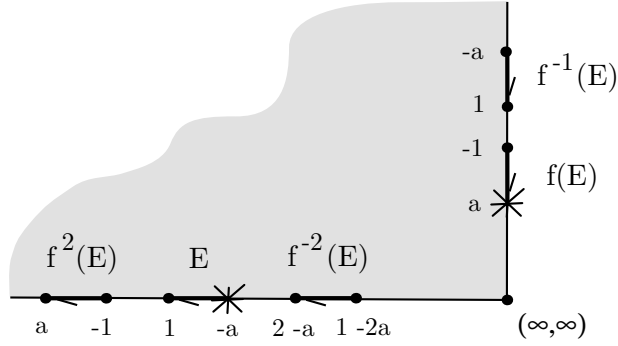}}
\caption{Orbit of the nontrivial rectangle
$R(\overline{10}\cdot\overline{10})=E$}
\end{figure}

In Figure 6 we have chosen one of the four ``quadrants'' abutting on
$(\infty,\infty)$ to illustrate the interval $E$ and part of its orbit.  The
``whisker'' coming off of $E$ indicates an orientation.  The point marked ``$*$''
in $E$ indeterminate for $f$; the other ``$*$'' is indeterminate for $f^{-1}$. 
Now recall the right hand side of Figure 4 in which $f(E)$ appears as a
vertical boundary segment of $R(00\cdot)$.  Figure 4 also shows some of
the canonical u-arcs that foliate $R(00\cdot)$.  The pictures of
$R((10)^k0\cdot)$ are similar except that as the number of digits
increases, the uniform u-arcs get closer to the limiting curve $\{y=-1\}$. 
The convergence is not uniform because the right endpoint of every
canonical u-arc is $(\infty,a)$, regardless of the number of digits. 
Hence, in the limit, the canonical u-arcs converge to the ``L''-shaped
rectangle $R(\overline{10}\,0\cdot)=f(E)\cup(R_0\cap\{y=-1\})$.  In light
of Corollary \ref{infpts}, it will follow from Theorem 6.5 that
$E=R(\overline{10}\cdot\overline{01})$.

\begin{thm}
\label{equivalentdef}
If $w\in\Sigma^*$, the two expressions on the right hand side of 
$$\dot R(w):=\bigcap_{k=-n}^m f^{-k} R_{w_k}=\bigcap_{k=-n}^m \dot f^{-k}
 R_{w_k}$$
are equal. 
Further, $\dot R(w)=R(w)$ if $w$ is alternating, and if $w$ is not alternating, we
have
$\dot R(w)=R(w)\cup\{(\infty,\infty)\}$ and $
R(w)=\dot R(w)-\{(\infty,\infty)\}$. In particular,
$R(w) = R(w^+)\cap R(w^-)$.
\end{thm}

\proof Let us start with the observation:
$$\dot f(R_0)-f(R_0)=\{y=-1\}-R_0\subset \interior R_+.$$
Thus 
$$(R_0\cup R_1)\cap\bigcap_{k=-n}^{-1}\dot f^{-k}R_{w_k}=
(R_0\cup R_1)\cap\bigcap_{k=-n}^{-1} f^{-k}R_{w_k},$$
from which we deduce that the definition of $\dot R(w)$ is unambiguous.  Next we
note that $(\infty,\infty)$ belongs to $R(w)$ if and only if $w$ is alternating. 
Finally, consider a point $p\in\overline{{\bf R}^2}-(\infty,\infty)$.  If the
orbit of $p$ is disjoint from the indeterminacy set, then $p\in R(w)$ if and only
if $f^jp\in R_{w_j}$ for all $j\in{\bf Z}$.  Thus $p\in R(w)$ if and only if
$p\in\dot R(w)$.  Otherwise, we may assume that $f^jp\in I(f)\cap (R_0\cup
R_1)=(-a,\infty)$ for some
$j\ge0$.  This case is handled by considering the various possibilities in
Corollary \ref{infpts}.
\qed

By the following result, $R$ is essentially a semi-conjugacy from
$(\sigma,\Sigma)$ to $(f,\cK)$.
\begin{thm}
If $w$ and $\sigma w$ are admissible, and if $(-a,\infty)\notin R(w)$ then
$R(\sigma w)=fR(w)$.
\end{thm}

We define 
$$\Omega:=\bigcup_{w\in\Sigma}
R(w).$$
\begin{thm}
$$\Omega=\bigcap_{n\in{\bf Z}} f^n(R_0\cup R_1).$$
\end{thm}

\proof The inclusion $\subset$ is evident.  We
will show the reverse containment.  For this, it suffices to show that
$$\bigcup_{w\in\Sigma}R(w)
\supset\bigcap_{m,n=0}^\infty\bigcup_{v\in\Sigma[-n,m]}R(v),$$
where $\Sigma[-n,m]$ denotes the set of admissible $[-n,m]$ words.   Now we
suppose that $p$ belongs to the right hand intersection.  Thus for
each $n,m$, there is a word $v$ of extent $[-n,m]$ with
$p\in R(v)$.  Let us suppose first that $p\notin I_+\cup I_-$.  If $n'<n''$ and
$m'<m''$, and if $v'$ and $v''$ are the corresponding words, then $v''$ extends
$v'$.  Thus there is a word
$w\in\Sigma$ of infinite length which is the common extension of all these finite
words.  It follows that $p\in R(w)$.  If $p\in I_+\cup I_-$, then by Corollary
\ref{infpts} we have $p\in R(w)$, where
$w$ is a finite subword of
$*0\overline{01}$.  In this case, too, we obtain an infinite word
$w\in\Sigma$ with $p\in R(w)$.  This gives the reverse containment, which
completes the proof. 
\qed

\section{Invariant cone fields; boundaries of rectangles}
\label{boundary}
\noindent  In this Section, we show the existence of invariant cone fields for
$f$.  This allows us to obtain slope bounds for s- and u-arcs.  From this we are
able to work more effectively with the boundaries of rectangles.

\begin{figure}
\centerline{\epsfxsize4.6in\epsfbox{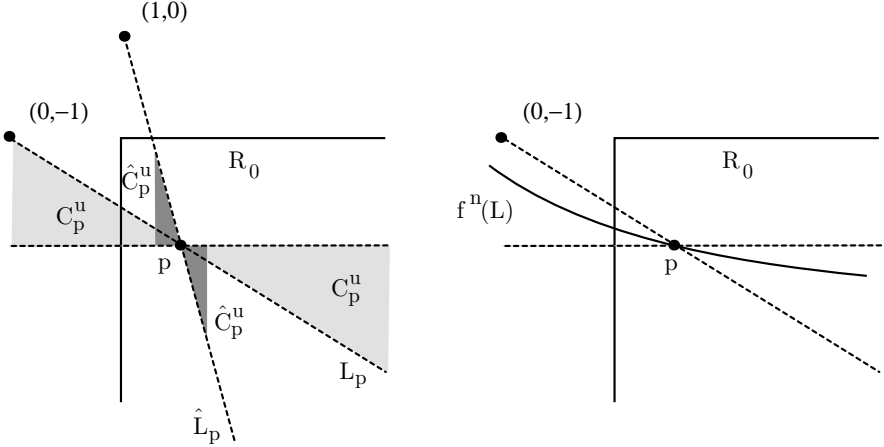}}
\caption{Tangent to $f^n(L)$ is between dashed lines.\label{slopebounds}}
\end{figure}

For a point  $p\in {\bf
R}^2$, we let $L_p$ denote the line from $(0,-1)$ to $p$ and $\hat L_p$ denote
the line from
$(1,0)$ to $p$.  Let $H_p$ and $V_p$ denote the horizontal and vertical lines
through $p$.  For
$p\in{\bf R}^2\cap R_0$, we let $\cC^u_p$ denote the cone of tangent vectors $t\in
T_p{\bf R}^2$ which are obtained by passing, in the counter-clockwise direction,
from $L_p$ to $H_p$.  In other words, $\cC^u_p$ contains those vectors in the
second and fourth quadrants between $L_p$ and $H_p$.  The cone
$\hat\cC^u_p$ is obtained by starting at $\hat L_p$ and passing in the
counter-clockwise direction until we reach $H_p$.   If
$q\in R_1$, then we let $\cC_q^u$ (respectively, $\hat\cC_q^u$) be the cone swept
out by starting at 
$V_q$ and moving counter-clockwise until we reach $L_q$ (respectively, $\hat
L_q$).  The cones $\cC^s$ and $\hat\cC^s$ are obtained as the images of $\cC^u$ and
$\hat\cC^u$ under the involution $(x,y)\mapsto(-y,-x)$.  Thus $\cC^s_p$
(respectively,
$\hat\cC_p^s$) is the complement of the interior of $\hat\cC^u_p$ (respectively,
$\cC_p^u$).  Figure \ref{slopebounds} shows
both cones for a point  $p\in R_0$; the corresponding picture for $p\in R_1$ is
obtained by reflecting about the line $y=x-1$.

\begin{thm}
If $p, fp\in (R_0\cup R_1)\cap{\bf R}^2$, then the differential
$Df_p$ maps vectors of
$\hat\cC_p^u$ to vectors in
$\cC_{fp}^u$.  
\end{thm}

\proof  Let us assume that $p\in R_0$ and fix a vector $t=(1,-\alpha)\in\hat\cC^u_p$
for
$\alpha>0$.  Let
$M=\{p+\zeta(1,-\alpha):\zeta\in{\bf C}\}$ denote the complex line passing through
$p$ in the direction $t=(1,-\alpha)$.  With respect to the basis
$\{\gamma_1^*,\gamma_2^*\}$ from \S3, the cohomology class $\{M\}\in H^2({\bf
P}^1\times{\bf P}^1)$ is the vector
$[1,1]$. Likewise, 
$\{fM\}=f_*\{M\}=[2,1]$.  We let $H_{fp}$ denote the horizontal line passing through
$fp$.  Then 
$\{H_{fp}\}=[1,0]$, and the intersection
multiplicity is
$$H_{fp}\cdot fM=[1,0]\cdot[2,1]=1.$$  
It follows that the intersection of $H_{fp}$
and $fM$ at $fp$ is transverse.  That is, 
$Df_p(t)$ is not horizontal.

Similarly, since $L_{fp}$ is neither vertical nor horizontal, we have
$\{L_{fp}\}=[1,1]$, and so
$$L_{fp}\cdot fM=[1,1]\cdot[2,1]=3.$$
Since $M\cap\{x=-a\}\ne\emptyset$, we have $(0,-1)\in fM$.  Thus each point of the
intersection
$$\{(0,-1),fp,(\infty,\infty)\}\subset L_{fp}\cap fM,$$ 
must have multiplicity one, which is to say that each intersection is transverse. 
Since
$\cC^u_{fp}$ is bounded by the horizontal and $L_{fp}$, we conclude that
$Df_pt\notin\partial\cC^u_{fp}$.

We may consider $fM$ as the union of arcs $\gamma_+:=f(M\cap\{x>1\})$ and
$\gamma_-:=f(M\cap\{x<1\})$.   As $t\to\pm\infty$, we have
$fp+t(1,-\alpha)=t(-\alpha,1)+O(1)$.  Thus $fM$ intersects $(\infty,\infty)$
through the second and fourth quadrants.  Further, since $t\in\hat\cC_p$, it follows
that $M\cap\{x=1\}\subset\{y\le0\}$.  Thus $\gamma_+$ begins at 
$(-\infty,+\infty)$ (in the second quadrant), passes through $(0,-1)$ and then
proceeds to
$(+\infty,a)$ (i.e.\ the
$y$-coordinate approaches $a$ as $x\to+\infty$).  We have seen that $L_{fp}$
intersects $\gamma_+$ transversally, and only in the points $(0,-1)$ and
$(\infty,\infty)$.  If
$fp$ is to the right of the point $(0,-1)$, then the portion of $\gamma_+$ to the
right of $fp$ must be above $L_{fp}$ because $\gamma_+$ approaches $(a,+\infty)$,
whereas $L_{fp}$ approaches $(+\infty,-\infty)$.  Thus the tangent to
$fM$ lies above the tangent to $L_{fp}$.  On the other hand, since the horizontal
$H_{fp}$ intersects $\gamma_+$ only once, the tangent to $fM$ at $fp$ must lie
below the horizontal.  The two other cases: $fp\in\gamma_+$ to the left of
$(0,-1)$, and $fp\in\gamma_-$ are handled similarly. \qed

\begin{thm}  If $p, fp, f^2p\in(R_0\cup R_1)\cap{\bf R}^2$,
then
$Df^2_p$ maps $\cC^u_p$ strictly inside $\cC^u_{fp}$.
\end{thm}

\proof   Since $\cC^u_p\subset\hat\cC^u_p$, it follows from the previous Theorem
that $\cC^u_p$ is mapped to $\cC^u_{fp}$.  Now we show that it is mapped strictly
inside.  Thus, if
$t\in\cC^u_p$, we must show that
$Df_pt\notin\partial\cC^u_{fp}$.  We may assume that $p\in R_0$; otherwise, we
work with $fp\in R_0$ instead.  In this case, $\partial\cC^u\cap\partial\hat\cC^u$
is horizontal, so it suffices to show that if $t$ is horizontal, then $Df_p(t)$ is
not in $\partial\cC^u_p$.  For this, let $H_p$ denote the horizontal
complex line passing through
$p$, with cohomology class $\{H_p\}=[1,0]$.  Thus
$\{fH_p\}=f_*\{H_p\}=[1,1]$, and so 
$$H_{fp}\cdot fH_p=[1,0]\cdot[1,1]=1$$
$$V_{fp}\cdot fH_p=[0,1]\cdot[1,1]=1,$$
where $V_{fp}$ denotes the vertical passing through $fp$.
It follows that the tangent to $fH_p$ is neither horizontal nor vertical at $fp$. 
Similarly, we have
$\{(0,-1),fp\}\subset L_{fp}\cap fH_p$, and the intersection number is 
$$L_{fp}\cdot fH_p=[1,1]\cdot[1,1]=2.$$  Thus
the intersection of $L_{fp}$ and $fH_p$ is transverse at $fp$, and so the tangent
does not belong to the boundary of $\cC^u_{fp}$.  We conclude that $Df_p(\cC^u_p)$
is strictly inside $\cC^u_{fp}$.
\qed

As a consequence, we obtain slope bounds on canonical s- and u-arcs; the first of
these slope bounds is illustrated in Figure \ref{slopebounds}.
\begin{thm}
\label{slope}
Let $L$ be a horizontal or vertical line which meets $R_0\cup R_1$ in a 
proper u-arc.  Let $p\in f^nL$ for some $n>0$, $m$ be the 
slope of $f^nL$ at 
$p$, and $\ti m$ be the slope of the line joining $p$ to $(0,-1)$. 
\begin{itemize}
\item If $p\in R_0$, then $\ti m < m < 0$.
\item If $p\in R_1$, then $m < \ti m < 0$.
\end{itemize} 
Suppose instead that $L$ meets $R_0\cup R_1$ in a proper s-arc
and that $\ti m$ is the slope of the line joining $p\in f^{-n}L$ to $(1,0)$.
\begin{itemize}
\item If $p\in R_0$, then $m < \ti m < 0$.
\item If $p\in R_1$, then $\ti m <  m < 0$.
\end{itemize}  
\end{thm}

We will call a u-arc $\gamma\subset R_0$, 
\emph{uniform} if it can be described as the graph $\{(x,g(x))\}$ 
of a function $g:[1,\infty] \to [-\infty,-1]$ such that 
 $g \equiv - \infty$, or
 $g$ is Lipschitz continuous with pointwise derivative
$g'(x)$ 
      constrained a.e.\ by the bounds in the first assertion of Theorem 
      \ref{slope}.
We extend the definition of uniformity to u-arcs in $R_1$ and s-arcs in $R_0$ 
and $R_1$ in the obvious fashion.  We say that an arc $\gamma\subset{\bf R}^2$ is
uniform if its closure is uniform.  Note, in connection with Figures 4 and 6, that
$R(\overline{10}\,0\cdot)\cap {\bf R}^2$ is a uniform u-arc.  The set
$R(\overline{10}\,0\cdot)$, however, is not a uniform u-arc, since it also contains
$f^{-1}(E)$.

With this terminology, we may summarize the 
first two assertions in Theorem \ref{slope} by saying that for any finite 
admissible $[-n,0]$ word $w$, the canonical u-arcs foliating $\interior R(w)$ 
are uniform.  The following is an easy consequence of the Arzela-Ascoli Theorem.

\begin{prop} 
\label{compact}  Let $\{\gamma_j\}_{j\in\N}\subset R_0$ be a sequence of 
uniform u-arcs with graphing functions $g_j$.  Suppose that 
$g = \lim_{j\to\infty} g_j$ exists pointwise on $[1,\infty)$. Then the 
convergence is uniform on compact subsets of $[0,\infty)$, and the limit 
extends to a function $g:[1,\infty]\to\R$ whose graph is a uniform u-arc in 
$R_0$. 
\end{prop}

\noindent Note that it is not necessarily the case that $g_j(\infty) \to
g(\infty)$.  This is  illustrated on the right hand side of Figure \ref{R00}, with
$g_j(\infty)=a$ and
$g(\infty)=-1$.

Uniformity of canonical arcs is the key to understanding the boundaries of
rectangles.

\begin{thm}
\label{dichotomy}
Let $w\in\Sigma_-$ be given.  Then for every $0\le n\le\infty$ the rectangle
$R(w[-n,0])$ is the set of points in $R_{w_0}$ between two uniform u-arcs
$\gamma_1$ and $\gamma_2$.  If $n<\infty$, then $\gamma_1\cap\gamma_2\cap{\bf
R}^2=\emptyset$.  If $n=\infty$, then either $\gamma_1=\gamma_2$, or $R(w)$ has
interior.
\end{thm}

\proof 
For the moment, suppose that $w=w[-n,0]$ is finite.  Then the canonical u-arcs 
$(\pi_w^u)^{-1}(t)$, $t\in T^u$, are uniform.  Therefore by Proposition 
\ref{compact} the two halves of $\partial^u R(w)$, which are pointwise limits 
of canonical u-arcs, are uniform u-arcs.  We pointed out earlier that these 
arcs meet, if at all, in a single infinite endpoint.

Now suppose that $w$ is infinite, and let $\gamma_{1,n},\gamma_{2,n}$ denote 
the uniform u-arcs bounding $R(w[-n,0])$.  Because $R(w[-n,0])$ decreases as 
$n$ increases, the graphing functions for $\gamma_{1,n},\gamma_{2,n}$ are 
monotone in $n$.  We apply Proposition \ref{compact} again as $n\to\infty$ to 
extract limiting uniform u-arcs.  The convergence of the graphing functions 
is uniform except at infinity, so we conclude that $R(w)\cap \R^2$ is the set 
of points in $R_{w_0}\cap\R^2$ between $\gamma_1$ and $\gamma_2$.  

Either $\gamma_1$ coincides with $\gamma_2$ or the corresponding graphing 
functions differ at some point.  In the first case,  
$R(w)\cap\R^2 = \gamma_1 = \gamma_2$.  In the second, continuity implies that
the graphing functions differ on an entire interval.  It follows that the 
region $R(w)$ has interior.
\qed

If $w\in\Sigma_-$, $n$, and $\gamma_1$, $\gamma_2$ are as in Theorem
\ref{dichotomy}, when we set
$$\partial^uR(w[-n,0])=\gamma_1\cup\gamma_2.$$
More generally, we may decompose $\partial R(w)$ into 
$$\partial^uR(w)=\partial^uR(w^-)\cap R(w),\ \ \ \partial^sR(w)=\partial^sR(w^+)\cap
R(w).$$
Since $\partial^uR(w)$ is a pair of (not necessarily distinct) uniform u-arcs, it
is natural to refer to the intersection of one of these arcs with $R(w)$ as a {\em
half} of $\partial^uR(w)$.
\begin{thm}
If $w$ is any admissible word, then
$$
\partial R(w) = \partial^u R(w) \cup \partial^s R(w).
$$
Each half of $\partial^u R(w)$ is connected and meets each half of 
$\partial^s R(w)$ in exactly one point.  If $\sigma w$ is well-defined and 
$R(w), R(\sigma w)\subset \R^2$, then
$$
f \partial^u R(w) = \partial^u R(\sigma w),\quad 
f \partial^s R(w) = \partial^s R(\sigma w).
$$
\end{thm}

\proof
The first conclusion is immediate from the fact that $R(w)=R(w^+)\cap R(w^-)$. 
Because of the bounds on slopes, a uniform u-arc in $R_{w_0}$ meets a uniform 
s-arc in $R_{w_0}$ in exactly one point.  Therefore, the second conclusion 
also follows.  Finally, if $R(w), R(\sigma w)\subset \R^2$, then $f$ maps 
$\partial R(w)$ homeomorphically onto $\partial R(\sigma w)$.  Moreover if $w$
is finite, the pair of proper u-arcs that make up $\partial^u R(w)$ must map 
to proper u-arcs.  Hence $f \partial^u R(w) \subset\partial^u R(\sigma w)$.  
Likewise, $f^{-1}\partial^s R(\sigma w)\subset \partial^sR(w)$.  This 
justifies the last conclusion for $w$ of finite extent.  A limiting argument 
justifies it for words of infinite extent.
\qed

We refer to the (at most four) points in $\partial^u R(w)\cap\partial^s R(w)$ 
as \emph{corners} of $R(w)$.  We denote the corner closest to the origin by 
$\delta R( w)$ and the corner 
furthest from the origin by $\ti \delta R( w)$.  Since uniform arcs of either type 
are graphs of non-increasing functions, we see that $\delta R( w)$ is also the 
corner nearest to the $x$-axis and to the $y$-axis and that $\ti\delta R( w)$ is 
likewise furthest from either axis.  Finally, if $w$ is a $[-n,m]$ 
word, we have $\delta R( w) \neq \ti\delta R( w)$ unless $n=m=\infty$.  When $n$ and
$m$ are both infinite, $\delta R( w) = \ti\delta R( w)$ if and only if $R(w)$ is a 
single point.

\section{Periodic points}
\label{periodic}
\noindent In this section we show that if $w\in\Sigma$ is periodic, then
$R(w)$ consists of a single periodic point (Theorem 6.3).  Further, the
correspondence $w\mapsto p\in R(w)$ is essentially a bijection between
periodic points of
$\sigma$ and periodic points of $f$.  Finally, by Theorem 6.4, all periodic
points except $(\infty,\infty)$ are of saddle type.

The alternating words $w\in\Sigma$ are special, as is the parabolic fixed
point $(\infty,\infty)$.  So we remove them from our discussion of
fixed points and define
$$ 
\text{Fix}'(\sigma^n)=\{w\in\Sigma:\sigma^nw=w\}-
  \{\overline{01},\overline{10}\}
$$ 
and
$$
\text{Fix}'(f^n)=\text{Fix}(f^n)-\{(\infty,\infty)\}=
\text{Fix}(f^n)\cap{\bf R}^2.
$$

\begin{prop}
\label{altfp}
If $w\in\Sigma$ is alternating, then $R(w)=\{(\infty,\infty)\}$.
\end{prop}

\proof Since $(x,y)\mapsto (-y,-x)$ conjugates $f$ to $f^{-1}$, the set 
$R(w) = R(w^+)\cap R(w^-)$ is symmetric about the line $y=-x$.  Therefore 
assuming $R(w)$ contains points other than $(\infty,\infty)$, we deduce that 
the corner $\delta R(w) \neq (\infty,\infty)$ of $R(w)$ opposite 
$(\infty,\infty)$ lies in $\R^2$.  By Theorem \ref{whatsinr} $R(w)$ avoids 
$I(f^n)$ for all $n\in\Z$.  Hence $f^2 R(w) = R(\sigma^2 w) = R(w)$ 
with corners sent to corners.  So because $f^2$ preserves $(\infty,\infty)$ 
it also preserves the opposite corner $\delta R(w) = f^2 \delta R(w)$.  
Lefschetz fixed point formula predicts that $f^2$ has five fixed points, and 
Proposition \ref{isolated} shows that $(\infty,\infty)$ accounts for four of 
these.  The point $((1-a)/2,(a-1)/2)\in R_0$ is fixed by 
$f$ (and therefore by $f^2$) and is thus the fifth fixed point.  Thus 
$\delta R(w) = ((1-a)/2,(a-1)/2)$.  But this cannot be, because by Theorem 
\ref{whatsinr}, $((1-a)/2,(a-1)/2)\in R(w)$ only for 
$w = \overline{0}$.  This contradiction shows that 
$R(w) = \{(\infty,\infty)\}$.
\qed

\begin{lem} 
\label{fp}
If $w\in\Sigma$ satisfies $\sigma^nw=w$, then $R(w)\cap{\bf R}^2$ contains
a fixed point for $f^n$.  This
point belongs to $\text{Fix}'(f^n)$ unless $w$ is alternating.
\end{lem}

\proof By Proposition \ref{altfp} we may suppose that $w$ is not alternating.  
Hence none of the subwords $v^- \eqdef w[-n,0]$, $v^+ \eqdef w[0,n]$, or 
$v \eqdef w[-n,n]$ is alternating.  Corollary \ref{infpts} therefore implies 
that $R(v)\subset \R^2$.  So we invoke Theorem \ref{interior2} to obtain that 
$\pi_{w_0}\circ \pi_v^{-1}$ maps $T^s_{w_0}\times T^u_{w_0}$ continuously and 
injectively into $T^s_{w_0}\times T^u_{w_0}$.  By Brouwer's Theorem, we obtain
a fixed point $q = \pi_{w_0}\circ \pi_v^{-1}(q) \in T^s_{w_0}\times T^u_{w_0}$.
The point $p = \pi_v^{-1}(q)\cap R(v) = \pi_{w_0}^{-1}(q)\cap R_{w_0}$ lies in
$\R^2$.  Breaking $\pi_v$ and $\pi_{w_0}$ into components and taking
advantage of the fact that $p\notin I(f^n)$, we see 
$$
\pi_{v^-}^u(p) = \pi_v^u(p) = \pi^u_{w_0}(p) = \pi^u_{v^+}(p), \quad 
\pi_{v^+}^s(p) = \pi_v^s(p) = \pi^s_{w_0}(p) = \pi^s_{v^-}(p).
$$
Therefore,
$$
\pi_{v^+}(p) = \pi_{v^-}(p) = \pi_{\sigma^n v^-}(f^n p) = \pi_{v^+}(f^n p).
$$
Since $p \in \R^2$, we conclude that $p = f^n(p)$.  By Theorem 
\ref{whatsinr} $p\in R(w).$ 
\qed

\begin{thm} 
\label{saddles}
If $w\in\Sigma$ satisfies $\sigma^nw=w$, then $R(w)=\{p\}$ is a single
point satisfying $f^np=p$.  If $w$ is not alternating, $p$ has multiplicity
one, and
$w$ and $p$ have the same period.  Finally, the map
$w\mapsto R(w)$ defines a bijection between $\text{Fix}'(\sigma^n)$ and
$\text{Fix}'(f^n)$.
\end{thm}

\proof  If $w$ is alternating, then $R(w)=\{(\infty,\infty)\}$ by Proposition 
\ref{altfp}.  Otherwise, $w\in\text{Fix}'(\sigma^n)$, and by Lemma \ref{fp} we 
may choose a point $p=p(w)\in \text{Fix}'(f^n)\subset\R^2$.  If $w$ and 
$\tilde w$ are distinct words in $\text{Fix}'(\sigma^n)$, then Theorem 4.7
implies
$p(w)\neq p(\ti w)$.  In particular 
$\#\text{Fix}'(\sigma^n) \leq \#\text{Fix}'(f^n)$.

From the discussion of symbolic dynamics in \S2, we have
$$
\#\text{Fix}'(\sigma^n)=\#\text{Fix}(\sigma^n)=F_{n+1}+F_{n-1}
$$
if $n$ is odd and
$$
\#\text{Fix}'(\sigma^n)=\#\text{Fix}(\sigma^n)-2 =F_{n+1}+F_{n-1}-2
$$
if $n$ is even.  We may also count the periodic points of $f$.  By Proposition 
\ref{isolated} and the equation preceding it, we have (ignoring multiplicity 
on the left hand sides)
$$
\#\text{Fix}'(f^n)\le F_{n+1}+F_{n-1}
$$
when $n$ is odd, and
$$
\#\text{Fix}'(f^n)\le F_{n+1}+F_{n-1}-2
$$
when $n$ is even.  In either case, it follows that 
$\#\text{Fix}'(\sigma^n)=\#\text{Fix}'(f^n)$, and the correspondence 
$w\to p(w)$ is bijective.  Further, since the count of fixed points of $f^n$ 
without multiplicity coincides with the count with multiplicity, we conclude 
that each element of $\text{Fix}'(f^n)$ has multiplicity one.

Next suppose that $R(w)$ contains two periodic points $f^np=p$, and 
$f^m\tilde p=\tilde p$.  Then $w \in \text{Fix}'(\sigma^{nm})$ and 
$p,\tilde p \in \text{Fix}'(f^{nm})\cap R(w)$, contradicting the previous 
paragraph.  Hence $R(w)$ contains at most one periodic point.

Recall that the period of $w$ is the smallest $n$ for which $\sigma^n w=w$.  
In particular, $p(w)\neq p(\sigma^k w)$ for $0\leq k < n$, so the period of 
$w$ divides the period of $p(w)$.  Lemma \ref{fp} implies that the reverse is
also true, so that $w$ and $p(w)$ have the same period. 

Now we wish to show that $R(w)$ is a point for $w\in \text{Fix}'(\sigma^n)$.  
We have $R(\sigma^k w)\subset\R^2$ for every $k\in\Z$.  Therefore, from
the  discussion at the end of Section \ref{boundary} $f^k:R(w)\to
R(\sigma^k w)$  is a corner preserving homeomorphism.  Since there are at
most four corners  for $R(w)$, each must be a periodic point for $f$. 
This implies that $R(w)$  has only one corner, which occurs if and only if
$R(w)$ is a point.
\qed

A consequence of the proof is that $(\infty,\infty)$ is a fixed point of
multiplicity 2 for odd iterates of $f$.
\begin{thm}
Every finite periodic point for $f$ is of saddle type.
\end{thm}

\proof  By Theorem 5.2,  every fixed point $p = f^n(p)$
is simple.  That is, no eigenvalue of $Df^n(p)$ can be one.  Nor can there
be an eigenvalue that is a $k$th root of unity, since that would mean that $p$
has multiplicity greater than one as a point of period $nk$ for some $k>1$.
Since $f$ preserves the area form $\zeta=dx\wedge dy/(y-x+1)$, whose 
singularities are disjoint from $R_0'\cup R_1'$, the product $\det Df^n(p)$ of 
the eigenvalues of $Df^n$ is exactly one.  We conclude that either $p$ is a 
saddle point, or that $Df^n$ is conjugate to an irrational rotation.
However, this latter conclusion is inconsistent with the fact that $f$
preserves the cone field $\cC^u_p$. \qed

We say that a word $w\in\Sigma_+$ is {\em eventually alternating} if there
exists 
$k\geq 0$ such that $w[k,\infty]$ is alternating. 

\begin{thm} 
\label{thm3}
If $w\in\Sigma_+$ is alternating, then $R(w)$ is one of the following
\begin{itemize}
\item $R(0\cdot\overline{10}) = \{x=\infty\}\cap R_0$; 
\item $R(1\cdot\overline{01}) = \{y=\infty\}\cap R_1$.   
\end{itemize}
If $w$ is eventually alternating and $k\geq 1$ is the minimum number for which 
$w[k,\infty]$ alternates, then $R(w)\cap{\bf R}^2$ is a uniform
s-arc inside $f^{-k+1}\{x=1\}\subset \cC(f^k)$.
\end{thm}

\proof   Let us start with the word $w=0\cdot \overline{10}$.  By
Corollary \ref{infpts}, $R(w)-{\bf R}^2=\{x=\infty\}\cap R_0$.  If
$R(w)\ne\{x=\infty\}\cap R_0$, then $\partial^sR(w)\cap{\bf R}^2$ also
contains a uniform s-arc $\gamma$.  It follows that $\gamma$ contains
$(\infty,\infty)$. since $R(w)\cap I(f^k)=\emptyset$ for $k\ge0$, we have
$f^2(R(w))=R(w)\cap R(010\cdot)$.  In particular,
$f^2\gamma\subset\gamma$.  By Lemma 3.2, the fixed points of $f^2$ are
isolated, so we can choose a neighborhood $U$ of $(\infty,\infty)$ so that
$f^2$ has no fixed points in $U-(\infty,\infty)$.  We may assume that
$f^2(\gamma\cap U)\subset \gamma\cap U$ (the case $f^{-2}(\gamma\cap
U)\subset\gamma\cap U$ is similar).  Since $f^2$ has no fixed points
except $(\infty,\infty)$, each $p\in \gamma\cap U$ satisfies
$f^{2m}p\to(\infty,\infty)$ as $m\to\infty$.  This, however, contradicts
Theorem 1.3.  Thus $R(w)=\{x=\infty\}\cap R_0$.

If $w$ is eventually alternating, then there is a $k>0$ such that
$(\sigma^kw)^+$ is alternating.  By Theorem \ref{whatsinr}, the image
$f^k(R(w)\cap {\bf R}^2)$ is contained in $R(0\cdot\overline{10})$ or
$R(1\cdot\overline{01})$, which are intervals at infinity.  The only
points in ${\bf R}^2$ that are sent to infinity by $f$ are those in
$\{x=1\}\subset\cC$.
\qed

We say that a word $w\in\Sigma_+$ is \emph{pre-periodic} if it is
pre-periodic for
$\sigma^+$. 

\begin{thm}
\label{perarc}
If $w\in\Sigma_+$ is pre-periodic but not alternating, then 
$R(w)\cap{\bf R}^2$ is a uniform s-arc.
\end{thm}

\proof By Corollary  \ref{dichotomy} it suffices to show that $\interior
R(w)=\emptyset$.  Theorem \ref{thm3} allows us to assume that $w$ is not
eventually alternating, i.e.\ the  (eventual) period $n$ of $w$ is larger
than two.  Replacing $w$ with 
$\sigma^j w$ for $j$ large enough, we can assume that $w[-n,\infty]$ is 
periodic.  Therefore neither $w[-n,0]$ nor $w[0,n]$ alternates, and from 
Corollary \ref{infpts} we have
$$
R(\sigma^{kn} w) \subset R(w[-n,n]) \subset \R^2
$$
for all $k\in\N$.  
Moreover, the invariant 2-form $\zeta=(y-x+1)^{-1}dx\wedge dy$ satisfies
$$
C^{-1}dx\wedge dy\le|\zeta|\le Cdx\wedge dy
$$
on $R(w[-n,n])$. Hence 
$$
\Area (\interior R(w)) \leq 
C \int_{\interior R(w)}|\zeta| = 
C \int_{\interior R(\sigma^{kn} w)}|\zeta| \leq 
C\Area R(\sigma^{kn}w).
$$
But $R(\sigma^{nk}w)$ decreases to $R(\tilde w)$, where $\tilde
w\in\Sigma$  is the periodic extension of $w[-n,n]$.  So by Monotone
Convergence and Theorem 
\ref{saddles}, it follows that $\Area R(\sigma^{nk}w)$ decreases to $\Area
R(\tilde w) = 0$ as $k\to\infty$.  We  conclude that $\Area(\interior R(w))
= 0$ and thus
$\interior R(w) =
\emptyset$.
\qed

\section{Uniform arcs and one-sided words}
\label{parabolic}
\noindent  
The following is one of the main
results of this paper.  

\begin{thm}
\label{stablearcs}
If $w\in\Sigma_+$ is not alternating, then $R(w)\cap{\bf R}^2$ is a uniform 
s-arc.
\end{thm}

When $R(w)-{\bf R}^2$ is a single point, the conclusion of Theorem 7.1
simplifies to the statement that $R(w)$ is itself a uniform s-arc.  The
only time when this does not happen is when the block `00' appears exactly
once in $w$, in which case $R(w)-{\bf R}^2=f^{-j}E$ for some $j$, by
Corollary \ref{infpts}. 

The rest of this section will be devoted to the proof of Theorem 7.1;
because of Theorem \ref{thm3} we will assume throughout that $w$ is not
pre-periodic.  By \S\ref{boundary} we know that it  is sufficient to show
that the area of
$R(w)$ is zero.  The invariant area  form is a useful tool, but it is
singular at infinity.   So we need to study  orbits that accumulate at
infinity, and for this we analyze the behavior near  the parabolic point.  We first characterize the itineraries of  points in $\Omega$
with unbounded forward orbits.

\begin{lem}
\label{compactorbits}
Let $w\in\Sigma_+$ and $K\subset R(w)\cap{\bf R}^2$ be a given compact
set.  Then the  forward orbit $\{f^n K\}_{n\geq 0}$ is unbounded if and
only if $w$ contains arbitrarily long alternating subwords.
\end{lem}

\proof 
By Theorem \ref{whatsinr} $f^nK \subset R(\sigma^n w)$.   
The lemma is therefore a consequence of Theorem \ref{thm3}.
\qed

Recall the 
invariant 2-form 
$$\zeta=(y-x+1)^{-1}dx\wedge dy.$$
\begin{thm}
Let $w\in\Sigma_+$ be given, and suppose that there is a number $M$ such
that any alternating subword of $w$ has length no greater than $M$.  Then 
$R(w)$ is a uniform s-arc.
\end{thm}

\proof 
Theorem \ref{perarc} allows us to assume that $w$ is not eventually periodic.  
If $n\ge M+1$, then neither $(\sigma^nw)^+$ nor $(\sigma^nw)^-$ is alternating,
and so $R(\sigma^nw)$ is a compact subset of ${\bf R}^2$.  By Lemma 
\ref{compactorbits} there is a compact set $S\subset {\bf R}^2$ which contains 
$R(\sigma^nw)$ for all $n\ge M+1$.  Let $C$ be a constant such that
$\zeta$ satisfies
$$
C^{-1}dx\wedge dy\le |\zeta|\le Cdx\wedge dy
$$
on $S\cap (R_0\cup R_1)$.  So 
$$
0 < \int_{R(\sigma^{M+1} w)}|\zeta| = \int_{R(\sigma^nw)}|\zeta|\le
C\Area R(\sigma^nw))\le C\Area S < \infty
$$
for all $n\ge M+1$.  

On the other hand, since $w$ is not eventually periodic, $(\sigma^n w)^+$ is 
different for every $n\in\N$.  Therefore the rectangles $R(\sigma^nw)$
are  mutually disjoint.  Since $\bigcup_{n=M+1}^\infty
R(\sigma^nw)\subset S$, we  must have  $\Area R(\sigma^nw) =0$.  It
follows that $R(\sigma^nw)$ has no  interior for any $n\in\N$.  By
Corollary \ref{dichotomy}, this is all we need  to know.
\qed
 
The final and most delicate part of Theorem \ref{stablearcs} is

\begin{thm}
\label{wandering}
Suppose that $w \in\Sigma_+$ contains arbitrarily long alternating
subwords. Then $R(w)$ is a uniform s-arc.
\end{thm}

\proof For 
$n\geq n_0$ sufficiently large, neither $(\sigma^nw)^-$ nor $(\sigma^nw)^+$ 
alternates.  Thus $R(\sigma^nw)\subset \R^2$.   For every $n\ge n_0$,
let 
$(x_n,y_n)=\delta R(\sigma^nw)$ and $(\ti x_n,\ti y_n)=\tilde\delta
R(\sigma^nw)$ be  the vertices of $R(\sigma^nw)$ which are closest and
farthest from the origin  in $\R^2$.  
\begin{lem}
\label{diam}
For every $n\ge n_0$, we have
$$
|x_n-\ti x_n||y_n-\ti y_n| \geq c\norm{(x_n,y_n)}.
$$
for $c=\int_{R(w)}|\zeta|$.
\end{lem}

\proof 
We estimate:
\begin{eqnarray*}
&  &    \int_{R(w)} |\zeta| = 
            \int_{f^nR(w)}| \zeta| =\int_{R(\sigma^nw)}|\zeta|\\
   & \leq & \int_{x_n}^{\ti x_n} \int_{\ti y_n}^{y_n} \frac{ dx\,dy}
            {\norm{(x_n,y_n)}}
            = \frac{|\ti x_n-x_n||y_n-\ti y_n|}{\norm{(x_n,y_n)}}.\\
\end{eqnarray*} 
The inequality follows for two reasons.  First, we replace
$R(\sigma^nw)$ by the euclidean rectangle with vertices $\delta
R(\sigma^nw)$ and $\tilde\delta R(\sigma^nw)$.  Then we estimate
$|\zeta|$ using the inequality $|y-x+1|\ge||(x_n,y_n)||$ on $R_0\cup
R_1$. 
\qed

\begin{lem}
\label{nearinfty}
There exists a large $M > 0$ such that if $\max\{x_n,-y_n\}>M$ and
$\max\{x_{n+2},-y_{n+2}\}>M$,  then 
$$f^2(x_n,y_n) = (x_{n+2},y_{n+2}) \text{ and } f^2(\ti x_n,\ti y_n) =
(\ti x_{n+2},\ti y_{n+2}).$$
\end{lem}

\proof  The coordinates of the rectangle $R(\sigma^nw)$ which are
nearest to the origin are 
$(x_n,y_n)$, so it follows that all points of the rectangle are at
distance at least $M$ from the origin.  Since $f^2$ acts by translation
on $\{x=\infty\}$ and $\{y=\infty\}$, $f^2$ is approximately a
translation on 
$\{\max\{|x|,|y|\}>M\}$.  Since $f^2$ maps the rectangle $R(\sigma^nw)$
to
$R(\sigma^{n+2}w)$, since
$Df^2$ preserves the cone fields $\cC^{s/u}$, and since the tangents
to the sides of the rectangles lie inside these cone fields, it follows
that the nearest and farthest vertices
$\delta$ and
$\tilde\delta$ are preserved.  \qed

For each $n\in\N$, we set $m_n = y_n/x_n$ and $\ti m_n = \ti y_n/\ti x_n$.

\begin{lem}
\label{slopeorder}
If $\ti M>0$ is large enough and $\min\{x_n,-y_n\} > \ti M$, then 
\begin{itemize}
\item $m_{n+2} < m_n < 0$ and $\ti m_{n+2} < \ti m_n < 0$;
\item $m_n < \ti m_n$ implies that $m_{n+2} < \ti m_{n+2}$;
\end{itemize}
\end{lem}

\proof
By the hypothesis, the points $(x_n,y_n)$, $(\tilde x_n,\tilde y_n)$,
$(x_{n+2},y_{n+2})$, and $(\tilde x_{n+2},\tilde y_{n+2})$ all belong to
$R_0$ and are near
$(\infty,\infty)$.  The points $(\ti x_n,\ti y_n)$, 
$(\ti x_{n+2},\ti y_{n+2})$ therefore have the same properties.

The `slope function' $m(x,y) = y/x$ is meromorphic on $\pp$ with a simple pole
along $\{x=\infty\}$ and a simple zero along $\{y=\infty\}$.  Since $f^2$ 
preserves both of these sets, $m\circ f^2$ also has a simple pole along 
$\{x=\infty\}$ and a simple zero along $\{y=\infty\}$.  The latter function has
further zeroes and poles along $\cC(f^2)$, but these avoid the point
$(\infty,\infty)$.  Hence 
$h(x,y) \eqdef (m\circ f^2)/m$ is holomorphic in a neighborhood of 
$(\infty,\infty)$.  

We introduce the change of variables 
$(x,y) \eqdef \varphi(s,t) \eqdef (1/s,1/t)$, setting 
$(s_n,t_n) = \varphi(x_n,y_n)$ and 
$(\ti s_n,\ti t_n) = \varphi(\ti x_n,\ti t_n)$.  Then $s_n > \ti s_n > 0$ and 
$t_n < \ti t_n < 0$, and both points $(s_n,t_n), (\ti s_n,\ti t_n)$ are near
the origin.  Moreover, $h\circ\varphi$ is holomorphic near $(0,0)$ and  
Proposition \ref{approx} tells us that
$$
h\circ\varphi(s,t) = 1 + 2s - 2t + O(\norm{(s,t)}^2).
$$
In particular, $h(x_n,y_n) = h\circ\varphi(s_n,t_n)$ and 
$h(\ti x_n,\ti y_n) = h\circ\varphi(\ti s_n,\ti t_n)$
both exceed one.  This implies the first assertion of the lemma.  Moreover,
$$
h(x_n,y_n) = h\circ\varphi(s_n,t_n) > h\circ\varphi(\ti s_n,\ti t_n) 
                                    = h(\ti x_n,\ti y_n),
$$
which implies the second assertion.
\qed

Now we complete the proof of Theorem \ref{wandering}.  Let
$c=\int_{R(w)}|\zeta|$ be the constant from Lemma 7.5. We will show
that $c=0$, so $\interior R(w)=0$.  The theorem will then follow from
Corollary 5.5.  So let us suppose, to the contrary, that $c>0$.  Choose
$j_0$ such that
$$j_0>\frac{|a|+1}{c}.$$
Let $\tilde M$ be as in Lemma 7.7, and increase $j_0$ if necessary to
obtain
$$R(0\cdot(10)^j)\subset\{\tilde M<x\}\text{ and
}R((01)^j\cdot0)\subset\{y<-\tilde M\}$$
for $j\ge j_0$.
The map $f^{2j_0}$ acts as translation by $j_0(a-1)$ on the line
$\{x=\infty\}$, so we may choose $M$ sufficiently large that the second
coordinate $\pi_2f^{2j_0}$ satisfies
$$|\pi_2f^{2j_0}(x,y)-y-j_0(a-1)|\le 1$$
for $(x,y)\in\{M<x,a\le y\le -1\}$.  Now that we have chosen $M$, we may
choose $k_0\ge j_0$ such that
$$R(0\cdot(10)^{k_0/2})\subset\{M<x\}\text{ and
}R((01)^{k_0/2}0\cdot)\subset\{y<-M\}.$$

The word $w$ contains arbitrarily long alternating subwords but is not
eventually alternating, so we may find $K\ge k_0$ and $N$ such that
$w[N,n+2K]$ is alternating and such that $w[N-1,N]=00$ and
$w[N+2K,N+2K+1]=00$.

For convenience of notation, let us suppose that $N=0$.  Thus $(x,y)\in
R(00\cdot)\subset\{1\le x<\infty,a\le y\le-1\}$ and $(x,y)\in
R(0\cdot(10)^K)\subset\{M<x\}$ since $K\ge k_0$.  By our estimate on the
second coordinate of $f^{2j_0}$, we have
$$a+j_0(a-1)-1\le y_{2j_0},\tilde y_{2j_0}\le -1 +j_0(a-1)=1.$$
By Lemma 7.5, we have
$$\tilde x_{2j_0}\ge x_{2j_0}+ \frac{c||(x_{2j_0},y_{2j_0})||}
{|y_{2j_0}-\tilde y_{2j_0}|} \ge x_{2j_0}+\frac{c}{|a|+1}x_{2j_0}.$$
Now we
estimate the slopes
$$m_{2j_0}=\frac {y_{2j_0}} {x_{2j_0}} \le 
\frac{2j_0(a-1) } {x_{2j_0}} $$
$$\tilde m_{2j_0}= \frac {\tilde y_{2j_0}} {\tilde x_{2j_0}}
\ge
\frac {a+j_0(1-a)-1} {\tilde x_{j_0}}  
>
\frac{ (j_0+1)(a-1) } {( 1+c/ (|a|+1))x_{2j_0}}.$$
By our choice of $j_0$, we have $m_{2j_0}<\tilde
m_{2j_0}\le 0$. Further, we have 
$$(x_{j},y_j)\in\{\tilde M<\min(x,-y)\}$$
for $2j_0\le j\le 2K-2j_0$.  Thus we may apply Lemma 7.7 to conclude that
$m_j<\tilde m_j$ for $2j_0\le j\le 2K-2j_0$.

On the other hand, we could have chosen the point $(x_{2K},y_{2K})$ as
our starting point.  In this case we use $f^{-1}$ instead of $f$ and
work backwards. (Passing from $f$ to
$f^{-1}$ corresponds to applying the involution
$(x,y)\mapsto(-y,-x)$.)  Our starting point satisfies
$(x_{2K},y_{2K})\in R(0\cdot0)\subset\{1\le x\le -a,y\le -1\}$ and
$(x_{2K},y_{2K})\in R((10)^K\cdot)\subset\{y\le-M\}$.   However, when we
perform the corresponding slope estimates, we obtain $m_{2K-2j_0}>\tilde
m_{2K-2j_0}$. 
 From this contradiction we conclude that $c=0$,
completing the proof of Theorem 7.4.  \qed

\section{Conjugacy with the subshift}
\label{conjugacy}
\noindent
We can now  make completely explicit  the connection between $f$ and the
golden mean subshift. The map $R$ turns out to be very nearly a topological
conjugacy.  Using this map we then transfer the unique measure of maximal entropy
from $\Sigma$ to $\Omega$ and draw a number of conclusions about the dynamics of
$f$.

Let $$\Sigma'= \Sigma-(W^s_{loc}(\overline{01},\overline{10})\cup
W^u_{loc}(\overline{01},\overline{10}))$$ 
denote the collection of those
words
$w$ such that neither $w^+$ nor $w^-$ is alternating. 

\begin{thm} 
\label{finiteomega}
For each $w\in\Sigma'$, the set $R(w)$ consists of a single point, and the assignment
$
w\mapsto R(w)
$
gives a homeomorphism between $\Sigma'$ and $\Omega\cap {\bf R}^2$.
\end{thm}

\proof
Since $w^+$ is not alternating, Theorem \ref{stablearcs} implies that
$R(w^+)$ is a uniform s-arc.  Similarly, $R(w^-)$ is a 
uniform u-arc.  Therefore, 
$$R(w)\cap {\bf R}^2 = R(w^+)\cap R(w^-)\cap{\bf R}^2$$ 
contains a unique 
point $p$.  Since $R(w)$ is connected (Corollary \ref{connected}), 
we conclude that $R(w)=R(w)\cap{\bf R}^2=\{p\}$.

If $w\in\Sigma$ and, for example, $w^+$ is alternating then 
$$R(w)\subset R(w^+)\subset \{x=\infty\}\cup\{y=\infty\}$$ by Theorem \ref{thm3}.
Therefore, we have from Corollary \ref{equivalentdef} that
$$
\Omega\cap{\bf R}^2 = \bigcup_{w\in\Sigma'} R(w).
$$
In other words, $w\mapsto R(w)$ is surjective.  But rectangles corresponding to
distinct admissible words intersect only at points outside $\R^2$ (Theorem 
\ref{whatsinr}).  So the assignment is also injective.  

To see that it is continuous, suppose that a sequence $\{w^j\}\subset\Sigma$ 
converges to $w\in\Sigma$.  Then for each $k$, there exists $j_0$ such that 
$j\geq j_0$ implies that $w^j[-k,k] = w[-k,k]$.  Hence 
$R(w^j)\subset R(w[-k,k])$.  Continuity then follows from the definition of 
$R(w)$ as the decreasing intersection of the sets $R(w[-k,k])$ as $k\to\infty$.
 
Continuity of the inverse map follows from continuity of $f$ away from its 
indeterminacy set and the fact that $p\in\R^2$ belongs to $R(w)$ if and only 
if $f^kp\in R_{w_k}$ for all $k\in\Z$. 
\qed
\begin{cor}
\label{density}
Saddle periodic points of $f$ are a dense subset of $\Omega\cap{\bf R}^2$, and 
$\Omega\cap{\bf R}^2$ is a totally disconnected and perfect subset of ${\bf R}^2$.
\end{cor}

\proof
It is well-known that $\text{Per}'(\sigma)$ is dense in $\Sigma'$ and that
$\Sigma = \overline{\Sigma'}$ is a perfect set.  Therefore the theorem 
follows directly from Theorems \ref{saddles} and \ref{finiteomega}.
\qed

Let us set
$$\Sigma'':=\Sigma-(W^s(\overline{01},\overline{10})\cup  
W^u(\overline{01},\overline{10}))=\bigcap_{n\in{\bf Z}}\sigma^n\Sigma'.$$ 

\begin{cor} The assignment $w\mapsto R(w)$ defines a topological conjugacy between
$(\sigma,\Sigma'')$ and $(f,\cD_f\cap \Omega)$.
\end{cor}
The following result shows that $f$ is topologically expansive on $\Omega\cap{\bf
R}^2$.  (Recall that $f^2f$ acts as a translation, and thus is not expansive, on
$\overline{{\bf R}^2}-{\bf R}^2$.)
\begin{thm}
\label{expansivitythm}
There is an $\eta>0$ such that if $p,q\in\Omega\cap{\bf R}^2$ are
distinct points, then
$\sup_{n\in{\bf Z}}dist(f^np,f^nq)>\eta$. 
\end{thm}

\proof Let $dist$ denote a distance function on $\overline{{\bf R}^2}$.  
Fix $\eta>0$ such that  $\eta<dist(R_1,R_0\cap\{1\le x\le -a\})$ and
$\eta<dist(R_1,R_0\cap\{-a\le y\le -1\})$.  By Theorem
\ref{finiteomega} there are $w,v\in\Sigma'$ such that $p=R(w)$ and $q=R(v)$.  If
$p\neq q$, we must have $v\ne w$.  Without loss of generality we may assume that
$w_0\ne v_0$ and thus $w_0=0$ and $v_0=1$.  If $dist(p,q)<\eta$, and if $p\in
R_0$, $q\in R_1$, then we must have $p\in R_0\cap\{-a<x\}$.  Now
$f(R_0\cap\{x>-a\})\cap R_0=\emptyset$.  It follows that
$fp\in R_1$.  And by Proposition 1.2, $fq\in R_1$.   Thus $w_1=1$ and $v_1=0$. 
Again, if
$dist(fp,fq)<\eta$, we must have $fq\in R_0\cap\{-a<x\}$.   Repeating the previous
observation, we conclude that $w^+$ and $v^+$ are alternating sequences.  A
similar argument applied to $f^{-1}$ shows that $w$ and $v$ are alternating. This
is a contradiction, since by Theorem 8.1, we have $w,v\in\Sigma'$.
\qed

Recall from \S2 the measure $\nu$ on $\Sigma$ of maximal entropy (equal to
$\log\phi$).  This measure puts no mass on the 2-cycle
$\{\overline{01},\overline{10}\}$.
Since $\nu$ is finite and all points of
$W^{s/u}(\overline{01},\overline{10})-\{\overline{01},\overline{10}\}$ are
wandering for the restriction of $\sigma$, it follows that $\nu$ puts no mass on
$W^{s/u}(\overline{01},\overline{10})$.  Thus $\Sigma''$ is a set of full measure
for $\nu$, and it follows that $\mu:=R_*\nu$ is a probability measure on
$\cD_f\cap\Omega$ which inherits the key properties of $\nu$:
\begin{cor} 
\label{mu}
The measure $\mu$ does not charge $I(f)$, is $f$-invariant and mixing and has
entropy 
$\log\phi$.  Further,
$$
\mu =\lim_{n\to\infty} \frac{1}{\# \text{Fix}'(f^n)} \sum_{p\in
\text{Fix}'(f^n)}
      \delta_p
=\lim_{n\to\infty} \frac{1}{\# \text{Fix}(f^n)} \sum_{p\in
\text{Fix}(f^n)}
      \delta_p,
$$
so $\mu$ reflects the asymptotic distribution of (saddle) periodic points of $f$.
\end{cor}

\begin{prop}
Let $\lambda$ be a probability measure on $\overline{{\bf R}^2}$ with the following
weak invariance property:  For each Borel set $E$  there are
sets
$E'\subset\dot f(E)$ and $E''\subset\dot f^{-1}(E)$ such that
$\lambda(E)=\lambda(E')=\lambda(E'')$.  Then $\lambda$ puts no mass on
$\overline{{\bf R}^2}-\cD_f$.
\end{prop}

\proof  Let us recall that
$$\overline{{\bf R}^2}-\cD_f=\bigcup_{n\in{\bf Z}} \dot f^n(I(f)\cup I(f^{-1})).$$
First, by the invariance property, $\lambda$ can put no mass on $I(f)$.  For if
$p\in I(f)$ has positive mass, then $\lambda(f^{-n}p)=\lambda(p)>0$ for $n\ge0$. 
Thus $\lambda$ would have infinite mass, since $f^{-n}p$ is disjoint from
$f^{-m}p$ if $n\ne m$.  Similarly, $\lambda$ puts no mass on
$\bigcup_{n\ge0}(f^{-n}I(f)\cup f^nI(f^{-1}))$.

Finally, consider a Borel subset $E\subset \dot f^Np$ for $p\in I(f)$.  Without
loss of generality $E$ is disjoint from $\bigcup_{j\ge0} f^jI(f^{-1})$.  Thus
$f^{-N}$ is smooth on $E$, and $\lambda(E)=\lambda(f^{-N}E)=\lambda\{p\}=0$.  
\qed

Any measure $\lambda$ on
$\overline{{\bf R}^2}$ which is $f$-invariant in the sense of the previous
Proposition will live on
$\cD_f$.  By Theorem 1.2, $\lambda$ can put no mass on $\interior(R_+\cup
R_-)$.  By Theorem 4.11, all of the mass of $\lambda$ is on $\Omega$, and
thus $\lambda$ is carried by $\Omega\cap\cD_f$.  Thus it will be of the form
$\lambda=R_*\eta$ for some
$\sigma$-invariant measure $\eta$ on $\Sigma$.  From the fact that $\nu$ is the
unique measure of maximal entropy on $\Sigma$, we obtain:
\begin{cor}
$\mu$ is the unique measure of entropy $\ge\log\phi$ on $\overline{{\bf R}^2}$.
\end{cor}

We say that a bi-infinite
sequence $\hat x=(x_n)_{n\in{\bf Z}}$ is an $f${\em -orbit} if $x_{n+1}\in\dot
fx_n$ for all $n\in{\bf Z}$.  Let $X$ be a compact subset of $\overline{{\bf
R}^2}$.  By
$\hat X_f$ we denote the space of $f$-orbits $\hat x$ such that $x_n\in X$ for all
$n\in{\bf Z}$.  This is a compact subspace of the infinite product space $X^{\bf
Z}$.  We let
$\hat f$ denote the shift map on
$\hat X_f$, which means that $\hat f\hat x=\hat y$, where
$\hat x=(x_n)$ and $\hat y=(y_n)$ are sequences with $y_n=x_{n+1}$.  It follows
that $\hat f$ is a homeomorphism of $\hat X_f$.

Let $\pi:\hat X_f\to X$ be the projection defined by $\pi\hat x=x_0$.  If
$x\in\cD_f$, then $x$ is contained in a unique $f$-orbit
$\iota(x):=(f^nx)_{n\in{\bf Z}}$.  In fact, $\pi:\pi^{-1}(\cD_f)\to\cD_f$ is a
homeomorphism, and its inverse is given by $\iota$.  We may use $\iota$ to push
$\mu$ up to an invariant measure
$\iota_*\mu$ on $\iota(\cD_f)\subset\widehat {{\bf R}^2}_f$.

\begin{prop}  If $\lambda$ is an invariant probability measure on $\hat X_f$ with 
$X=\overline{{\bf R}^2}$, then $\lambda$ puts full measure on
$\iota(\cD_f)$.  Thus $\iota_*\mu$ is the unique measure of entropy $\log\phi$ on
$\iota(\cD_f)$.
\end{prop}
\proof  Pushing $\lambda$ down to $\overline{{\bf R}^2}$, we obtain a measure 
$\pi_*\lambda$ which is invariant in the sense of Proposition 8.5.  Thus
$\pi_*\lambda$ puts no mass on the complement of $\cD_f$.  Thus
$\lambda$ can put no mass on $\pi^{-1}(\overline{{\bf R}^2}-\cD_f)$.

If $\lambda$ is an invariant measure of entropy $\log\phi$ on
$\hat X_f$, then $\lambda$ lives on $\iota\cD_f$.  Thus
we may identify $\lambda$ with an $f$-invariant measure on $\cD_f$ with entropy
$\log\phi$.  This measure must be $\mu$, so $\lambda=\iota\mu$
\qed

Now we discuss the topological entropy of $f$.  The approach we follow here is to
replace $f$ by the map $\hat f$ acting on the orbit space.  In this case,  $(\hat
f,\hat\Omega)$ is a compactification of the restriction of $f$ to $\cD_f$.  A
second approach would be to work directly with
$f$, as is done by Guedj [Gu2].
\begin{thm}  The topological entropy of $\hat f$ on $\hat X_f$ is equal to
$\log\phi$ for $X=\Omega$, $\overline{{\bf R}^2}$, and ${\bf P}^1\times{\bf
P}^1$.  
\end{thm}

\proof  Let us consider first the case $X={\bf P}^1\times{\bf P}^1$.  In this case
Friedland [Fr] has shown that $h_{top}(f,\hat X_f)$ is bounded above by the
logarithm of the spectral radius of $f^*$ action on cohomology $H^*(X)$.  We have
seen that $f^*$ is represented by the matrix $
\left(\begin{matrix}
1&1 \\
1&0
\end{matrix}\right)
$, and thus the spectral radius is
given by the golden mean $\phi$.  It follows that $h_{top}(f,\hat X_f)\le\log\phi$
for all three choices $X$.

Now we consider the case $X=\Omega$.  We have seen that $\iota_*\mu$ is an
invariant measure on $\hat X_f$ with entropy equal to $\log\phi$.   Since the
topological entropy dominates the entropy of any invariant measure, it follows
that $h_{top}(f,\hat X_f)\ge\log\phi$.  Thus $h_{top}(f,\hat X_f)=\log\phi$ for all
three choices of $X$. \qed

We note that $(\infty,\infty)$ belongs to $\cD_f$ and is contained in a unique
orbit
$\iota(\infty,\infty)$, which is the constant sequence
$(\infty,\infty)$.  Let us write
$$\hat\Omega_*:=\hat\Omega_f-\iota(\infty,\infty)
=\pi^{-1}(\Omega-(\infty,\infty)).$$  
Each $x\in\Omega-(\infty,\infty)$ is contained in $R_j$ for a unique $j$.  Thus we
have a \emph{coding map}
$$c:\hat\Omega_*\to\Sigma$$ 
given by $c(\hat x)=(w_n)$, where $w_n$ is chosen such that $x_n\in R_{w_n}$ for
all $n\in{\bf Z}$.  It follows that
$$c:(\hat f,\hat\Omega_*)\to(\sigma,\Sigma)$$
is a semi-conjugacy.  This is an inverse to the mapping $R$ in the following
sense:  if $x\in\cD_f$ and $\hat x=\iota(x)$, then 
$R(c(\hat x))=\{x\}$, which means that 
$\pi=R\circ c$ as a mapping from $\iota(\cD_f)$ to $\cD_f.$
In other words,
\begin{prop} For $\hat x\in\pi^{-1}(\Omega\cap{\bf R}^2)$, we have $\pi(\hat
x)=R(c(\hat x))$.  More generally, for $\hat x\in\hat\Omega_*$, we have $\pi(\hat
x)\in R(c(\hat x))$.  
\end{prop}

\section{Parabolic basin; Nonwandering set}
\label{parabolicbasin}
\noindent   In this section and the next, we study the sets
$$\Omega_+=\bigcup_{w\in\Sigma_+}R(w), \ \text{ and }
\ \ \Omega_-=\bigcup_{w\in\Sigma_-}R(w).$$  By Theorem 4.9 we have
$$\Omega=\Omega_+\cap\Omega_-.$$ 
 By Theorem 7.1, each $R(w)$ is a uniform arc, and because of this, the sets
$\Omega_\pm\cap{\bf R}^2$ have product structure.  In this section we show that
$\Omega_+$ is the complement of the forward basin $\cB_+$ inside $R_0\cup R_1$.  Then we
identify the nonwandering set as the complement of the total basin $\cB_+\cup\cB_-$ inside
$\overline{{\bf R}^2}$. 

We define the {\em forward basin} of $(\infty,\infty)$, written $\cB_+$, to be the set of
points $p$ which are contained in neighborhoods $U$ such that $f^n|U$ converges uniformly to
$(\infty,\infty)$ as $n\to+\infty$. 

\begin{thm} The forward basin is
$\cB_+=int (W^s(\infty,\infty))$.
\end{thm}

\proof  Let us start by noting that $I_+\cap\cB_+=\emptyset$.  For if $p\in I(f)$, then
$\dot fp$ is one of the horizontal curves in the right hand side of Figure 1, and for all
$j\ge1$
$\dot f^jp$ is a nontrivial element of $\pi_1(\overline{{\bf R}^2})$.  Thus $\dot f^jp$ cannot
be contained in a disk about $(\infty,\infty)$.

It follows that $\cB_+$ is an open subset of $W^s(\infty,\infty)$, so it suffices to show
that $int(W^s(\infty,\infty))=\cB_+$.  In fact we will show that these sets both coincide
with the set of all points of $W^s(\infty,\infty)$ whose forward orbits are contained in
${\bf R}^2$.  One direction is clear.  For if $p\in
W^s(\infty,\infty)$, we must have
$f^np\in R_+$ for some positive $n$.  If $f^np\in R_+\cap{\bf R}^2$, then there is a
neighborhood
$U$ of $p$ such that $f^{n+1}U\subset R_+$.  Thus $p$ belongs to both $\interior
W^s(\infty,\infty)$ and
$\cB_+$.

It remains to consider the case of $p\in W^s(\infty,\infty)$ such that $f^np\notin{\bf R}^2$
for some
$n\ge0$.  Let us first claim that there exists a neighborhood $U$ of
$p$ such that $f^kU\cap R_0$ contains a connected open set $V$ such that $\overline
V\cap\{x=\infty\}\ne\emptyset$.  Let $n\ge0$ be the smallest number such that
$f^np\notin{\bf R}^2$.  If $n=0$, then for sufficiently large $k$, $f^kU$ will in fact
contain a one-sided neighborhood of $f^kp$ inside $R_0$.  If $n>0$, then
$f^{n-1}p\in\{x=1\}\cap{\bf R}^2$, and $f^{n-1}U$ is a neighborhood of $f^{n-1}p$ in ${\bf
R}^2$.  We see that $f^{n+k}U$ will have the desired property.

We will conclude the proof by showing that $V\not\subset W^s(\infty,\infty)$ and thus
$p$ belongs to neither $\interior W^s(\infty,\infty)$ nor $\cB_+$.  By Theorem 6.5,
$R(\cdot(01)^j)$ is a small one-sided neighborhood of $\{x=\infty\}\cap R_0$ inside $R_0$,
and this neighborhood shrinks to $\{x=\infty\}\cap R_0$ as $j\to\infty$.  Thus we may choose
$j$ sufficiently large that any s-arc in $R(\cdot(01)^j)$ crosses $V$.  On the other
hand, if $w$ is not eventually alternating, we have $R(w)\cap
W^s(\infty,\infty)=\emptyset$.  Since words which are not eventually alternating can begin
with
$(01)^j$, we see that $V\not\subset W^s(\infty,\infty)$, as claimed.  \qed

To define the {\em backward basin} $\cB_-$, we replace $f$ with $f^{-1}$; the analogue of
the theorem above holds for
$W^u(\infty,\infty)$ and $\cB_-$.

Each of the sets $\Omega_+\cap R_j\cap{\bf R}^2$ carries a natural product structure.  To
see this, let us define
$$\cT_{j,\pm}=\{w\in\Sigma_\pm:w_0=j,\text{ and } w \text{ not alternating}\}.$$
For instance, if $j=0$,  we set $I=(-\infty,-1]$, and we have a
homeomorphism
$$\cT_{+,0}\times I\ni(w^+,t)\mapsto R(w^+)\cap\{y=t\}\in\Omega_+\cap R_0.$$
To see that this map in fact defines a homeomorphism, note that for each 
$w^+\in\cT_{+,0}$, $R(w^+)\cap{\bf R}^2$ is a uniform s-arc in $R_0$ which
intersects $\{y=t\}$ in a single point.  Similarly,  if $j=1$, we set $I=[1,\infty)$ and 
have a homeomorphism $$\cT_{+,1}\times I\cong \Omega_+\cap R_1\cap {\bf R}^2.$$

\begin{lem}
$\Omega_+=(R_0\cup R_1)-\cB_+$.
\end{lem}

\proof  If $p\notin\Omega_+$, then by Proposition 1.2, $f^jp\in R_+-(R_0\cup R_1)$ for some
$j\ge1$.  It follows that the complement of $\Omega_+$ is contained in
$W^s(\infty,\infty)$.  If $w\in\Sigma_+$ is not eventually alternating, then $R(w)$ is
disjoint from $W^s(\infty,\infty)$. Since the set of all such $R(w)$ is a dense subset of
$\Omega_-$ and since
$\Omega_+$ is nowhere dense in $R_0\cup R_1$, the Lemma follows.  \qed

\begin{thm} $\overline{{\bf
R}^2}-\cB_+=\overline{{\bf R}^2}\cap\partial\cB_+=\bigcup_{n\ge0}f^{-n}\Omega_+$, and
$\cB_+$ is dense in ${\bf R}^2$.
\end{thm}

\proof 
By Proposition 1.2, $f^{-1}\Omega_+-\Omega_+\subset R_-$.  Thus
$f^{-n-1}\Omega_+-f^{-n}\Omega_+$ tends uniformly to $(\infty,\infty)$ through $R_-$ as
$n\to\infty$. The Theorem then follows from the previous Lemma.
\qed

A point $p$ is said to be {\em wandering} if it has a neighborhood $U$ such that $U\cap
f^n(U-I(f^n))=\emptyset$ for all $n\ne0$.
\begin{prop}
The points of $\cB_+\cup\cB_-$ are wandering.
\end{prop}

\proof  If $p\in\cB_+$, let $U$ be a neighborhood of $p$ on which $f^n|U$ converges uniformly
to $(\infty,\infty)$.  It follows that, shrinking $U$ if necessary, we have $f^nU\cap
U=\emptyset$ for $n\ge1$.  Since $p\notin I(f^n)$ for $n\ge0$, we also have that $f^nU\cap
U=\emptyset$ for all $n\le-1$.

Now suppose $p=(1,0)$.  We have $\dot fp=\{y=a\}\subset R_+\cup R_0$ (see Figures 1 and 2). 
By Proposition 1.2, $\dot f^np\subset R_+\cup R_0\cup R_1$ for all $n\ge1$.  It follows that
if
$U$ is a sufficiently small neighborhood of $p$, we have $f^n(U-\{p\})\subset R_+\cup
R_0\cup R_1$ for all $n\ge1$.  In particular, we may choose $U$, so that $U\cap
f^n(U-\{p\})=\emptyset$ for $n\ge1$.  On the other hand, since $(1,0)\in R_-$, we may choose
$U$ small and apply Proposition 1.2 to have $f^{-n}U\cap U=\emptyset$ for $n\ge1$.  \qed

\begin{prop}
The lines at infinity $\{x=\infty\}\cup\{y=\infty\}$ are contained in the
nonwandering set of $f$.
\end{prop}

\proof  Since the nonwandering set is closed and invariant, it suffices to
show that any point $p\in R_0\cap\{x=\infty\}-I_-$ is nonwandering.  Let
$U$ be a neighborhood of $p$.  Let us choose $N$ such that $f^{-N}p\in
R_0\cap\{y=\infty\}$.  By Theorem \ref{thm3},
$U$ intersects
$f^{-m}\{x=1\}$ for sufficiently large $m$.  Hence $f^{2m}U\cap R_0$
contains a component connecting $\{x=1\}$ to $\{x=\infty\}$.  Further,
since this component is inside $R((10)^{2m}\cdot)$ which is a small
one-sided neighborhood of $R_0\cap\{y=\infty\}$, it must cross $f^{-N}U$.
\qed

\begin{thm} The nonwandering set is 
$$\{x=\infty\}\cup\{y=\infty\}\cup\Omega =
\overline{{\bf
R}^2}-\left(\cB_+\cup\cB_-\right).$$
\end{thm}

\proof  
The two sets above are equal by Theorem 9.3.  By Proposition 9.4, this set
contains the nonwandering set.  The Theorem now follows from Proposition 9.5 and the
fact that the periodic (saddle) points are dense in $\Omega\cap{\bf R}^2$.
\qed

\section{Stable manifolds and laminar currents}
\label{currents}
\noindent This section is devoted to identifying the stable manifolds of points of
$\Omega$.  In order to do this, we define a pseudometric on $R_0\cup R_1$ which is
uniformly contracted/expanded by $f$.  Finally, we use the stable manifolds and
transverse measures to construct stable and unstable currents, as was done in the
Axiom A case by Ruelle and Sullivan [RS].  One motivation for doing this is to
present all of the aspects of hyperbolicity present in $f$.  Another motivation is
that laminar currents have proved increasingly useful in understanding the dynamics
of mappings (see, for instance [Du]).

We say that $p\in\Omega_+$ is \emph{accessible} if there exists a continuous curve
$\gamma:[0,1]\to(R_0\cup R_1\cup R_+)$ such that
$\gamma([0,1))\cap\Omega_+=\emptyset$, and $\gamma(1)=p$.

\begin{thm}
The accessible points of $\Omega_+$ are the points in the set
$$W^s(\infty,\infty)\cap\Omega_+=\Omega_+\cap\bigcup_{n\ge0}f^{-n}\{x=1\}$$
$$=(R_0\cup R_1)\cap \bigcup_{n\ge0}f^{-n}\{x=1\}.$$
\end{thm}

\proof 
Suppose that $p\in\Omega_+$ is accessible.  Without loss of generality we may
suppose that $p\in R_0$, and $\gamma\subset R_0$ is as in the definition.  Let us
assume, by way of contradiction, that $f^np\notin\{x=1\}$ for all $n\ge0$.  Since
$f^n(\gamma-\{p\})\cap\Omega_+$, it follows that
$f^n(\gamma)\cap(\{x=1\}\cap R_0)=\emptyset$.  We may choose $n$
sufficiently large that $f^np\in R_0$ and $f^n(\gamma)\cap R_+\ne\emptyset$, for
otherwise
$\gamma\subset\Omega_+$.  Since $f^n\gamma\subset R_0\cup R_1\cup R_+$ it follows
that $f^n\gamma$ contains a point of $\{x=1\}\cap  R_0$.  This contradiction 
shows that the accessible points are contained in the displayed set.

Conversely, suppose that $f^np\in\{x=1\}\cap R_0$.  Then $f^{n+1}p\in
R_0\cap\{x=\infty\}$, and so $p\in\Omega_+$.  Choose $\gamma$ to be a curve which
has $f^np$ as an endpoint and such that $\gamma-\{f^np\}\subset \interior(R_+)$. 
Then $f^{-n}(\gamma-\{p\})\cap\Omega_+=\emptyset$, and $f^{-n}\gamma$ is an access
to $p$.
\qed

Recall the homeomorphism
$\cT_{\pm,j}\times I\cong\Omega_{\pm,j}$ from the previous section.  With this, we
may define an ordering on
$\cT_{\pm,j}$.  For instance, if $w^1,w^2\in\cT_{+,0}$ (resp.\
$w^1,w^2\in\cT_{+,1}$) we say that $w^1<w^2$ if the s-arc $R(w^1)\cap{\bf R}^2$
lies to the left of (resp.\ below) $R(w^2)\cap{\bf R}^2$.  It is evident that this
defines total orders on $\cT_{+,0}$ and $\cT_{+,1}$.  If $w^1,w^2\in\cT_{+,j}$ and
$w^1\le w^2$, we define the interval $[w^1,w^2]=\{x\in\cT_{+,j}:w^1\le x\le w^2\}$.
By the product structure, we see that $[w^1,w^2]$ is a closed subset of
$\cT_{\pm,j}$, and $(w^1,w^2)$ is open.

\begin{figure}
\centerline{\epsfxsize4.6in
            \epsfbox{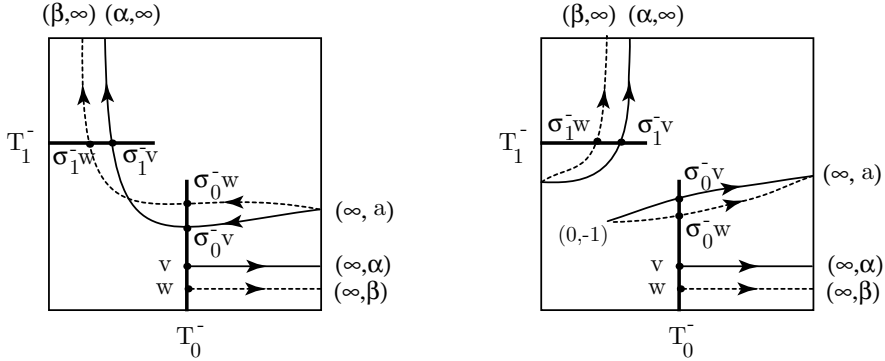}}
\caption{Induced mappings $\sigma^-_0$ and $\sigma^-_1$ on the transversal
${\cT}_-$.
\label{intervals}}
\end{figure}

The map $f^{-1}$ induces maps $\sigma^+_k:\cT_{+,0}\to\cT_{+,k}$, $k=0,1$ and 
$\sigma^+_0:\cT_{+,1}\to\cT_{+,0}$ as follows.  If $w=w_0\cdot
w_1w_2\dots\in\cT_{+,0}$, then for
$k=0,1$, $f^{-1}R(w)\cap R_k$ is the s-arc $R(k\cdot w_0w_1\dots)$.  We write
$\sigma^+_k(w)=k\cdot w_0w_1w_2\dots$; thus $\sigma^+_0$ and
$\sigma^+_1$ are the two right inverses of $\sigma_+$.
We define $\sigma^-_k$ in an analogous fashion.

Figure 8 shows how $f$ maps u-arcs in $R_0$.  From this we can see that the
restriction of
$\sigma_j^-$ to each $\cT_{-,k}$ is either order-preserving or order-reversing.
Thus we have:
\begin{prop}
The mappings $\sigma^\pm_k$ take intervals to intervals.
\end{prop}

The sets $\cT_{\pm,j}$ have the topologies of closed, perfect, totally
disconnected subsets of ${\bf R}$.  We may write the complement of $\cT_{\pm,j}$
in ${\bf R}$ as a union of open intervals $(a,b)$: the points
$w^a,w^b\in\cT_{\pm,j}$ corresponding to $a$ and $b$ are said to be
\emph{accessible} points of $\cT_{\pm,j}$.  In terms of the order structure on
$\cT_{\pm,j}$, these pairs of accessible points are characterized by the condition
that the intervals $[w^a,w^b]$ contain only the endpoints $w^a$ and $w^b$.

\begin{prop}
A point $w\in\cT_{+,j}$ is accessible if and only if it lies in the stable
manifold
$W^s(\{\overline{01},\overline{10}\},\sigma_+,\Sigma_+)$.  
\end{prop}

Now we define a distance function $\dist'_\pm$ on each of the sets $\cT_{\pm,j}$. 
For convenience of notation, let us work with $\cT_{+,0}$.  Let $\nu_{+}$ denote
the measure of maximal entropy for
$(\sigma_+,\Sigma_+)$.  For $w^1,w^2\in\cT_{+,0}$ we define
$$\dist'_+(w^1,w^2)=\nu_{+}([w^1,w^2])\ge0$$
If $w^1$ and $w^2$ are not both accessible, then $[w^1,w^2]$ contains an open
subset of $\Sigma_+$, so $\dist'(w^1,w^2)=\nu_+([w^1,w^2])>0$. An immediate
consequence of the balanced property of $\nu_+$ is:
\begin{prop}
For $v,w\in\cT_{\pm,k}$, 
$\dist'_\pm(\sigma^\pm_jv,\sigma^\pm_jw)=\phi^{-1}\dist'_\pm(v,w)$.
\end{prop}

Now we may extend $\dist'_\pm$ to $R_j$ by the following requirement: if $U$ is a
connected component of $R_j-\Omega_+$, then we set $\dist'_+(p,q)=0$ for all
$p,q\in U$.  Thus
$\dist':=\max(\dist'_+,\dist'_-)$ is a pseudodistance on $R_j$ for $j=0,1$.  Since
$R_0\cap R_1=\{(\infty,\infty)\}$, we extend $\dist'$ to $R_0\cup R_1$ by setting
$\dist'(p,q)=\dist'(p,(\infty,\infty))+\dist'(q,(\infty,\infty))$ whenever $p\in
R_0$ and $q\in R_1$.

Recall from Proposition 1.2 that if $p\notin W^s(\infty,\infty)\cup I_+$, then
$f^np\in R_0\cup R_1$ for $n$ sufficiently large.  This motivates
the following definition.  For
$p\in\overline{{\bf R}^2}-(W^s(\infty,\infty)\cup I_+)$, we set
$$W^s(p):=\{q\in\overline{{\bf
R}^2}-(W^s(\infty,\infty)\cup I_+):\lim_{n\to+\infty}\dist'(f^np,f^nq)=0\}.$$ 

\begin{thm}
Let $p\in\Omega-W^s(\infty,\infty)$ be given, and let $w\in\Sigma$ be such that
$p=R(w)$.  Then $\dot
W^s(p):=\bigcup_{n\ge0}f^{-n}R((\sigma^nw)^+)$ is equal to $W^s(p)\cup I_+$, and 
$W^s(p)=\dot W^s(p)-I_+\subset{\bf R}^2$. 
\end{thm}

\proof  On $R(w^+)$, $\dist'$ is equal to $\dist'_-$ and there decreases by
factors of $\phi$ under forward iteration, so 
$R(w^+)\subset W^s(p)$.  Thus $\dot W^s(p)-I_+\subset W^s(p)$.  

Conversely, let us
suppose that $q\in W^s(p)-(W^s(\infty,\infty)\cup I_+)$.  
Let $U$ be a neighborhood of
$(\infty,\infty)$ such that $f$ is a diffeomorphism from $U\cap R_0$ (resp.\
$U\cap R_1$) to
$fU\cap R_1$ (resp.\ $fU\cap R_0$).  Choose
$\epsilon>0$ such that an
$\epsilon$ neighborhood of $(\infty,\infty)$ with respect to $\text{dist}'$ is
contained in
$U$.  For $N$ sufficiently large, we may assume that
$\text{dist}'(f^np,f^nq)<\epsilon/2$ for all $n\ge N$.   By Proposition 1.2,
$f^np,f^nq\in R_0\cup R_1$ for $n$ sufficiently large.
We claim that $f^np$ and
$f^nq$ must be contained in the same rectangle for $n\ge N+1$.  Thus it will
follow that  $f^{N+1}q\in R(\sigma^{N+1}w)$, which means that
$q\in \dot W^s(p)$.  To see why the claim is true, observe that if $f^np\in
R_0\cap\Omega$ and $f^nq\in R_1\cap\Omega$ are points in different rectangles,
then since $f^np,f^nq\in U$ for all $n\ge N$, we have $f^{n+1}p\in R_1$ and
$f^{n+1}q\in R_0$.  Thus $f^np$ and $f^nq$ are in opposite $R_j$'s for $n>N$,
which means that $p$ and $q$ would have to corresponding to eventually alternating
words, which is not possible since $p\notin W^s(\infty,\infty)$.
\qed

\begin{cor}
If $p\in\Omega-W^s(\infty,\infty)$, then $W^s(p)\cap(R_0\cup
R_1)\subset\Omega_+$.  If $p\in\Omega-W^s(\infty,\infty)$ and
$q\in\Omega-W^u(\infty,\infty)$, then $W^s(p)\cap W^u(q)\subset\Omega$.
\end{cor}

Let us define 
$$\cW^s={\bf R}^2\cap \bigcup_{n\ge0}f^{-n}\Omega_+-I_+, \ \
\cW^u={\bf R}^2\cap \bigcup_{n\ge0}f^{n}\Omega_--I_-.$$
\begin{thm}
$\cW^s$ is a lamination of ${{\bf
R}^2}-I_+$; the leaves of $\cW^s$ are the connected components of
$W^s(p)\cap{\bf R}^2$ for $p\in\cD\cap\Omega$ and of ${\bf
R}^2\cap f^{-n}\{x=1\}-\cI_+$ for
$n\ge0$.  Similarly, $\cW^u$ is a lamination of
${{\bf R}^2}-I_-$; the leaves of $\cW^u$ are the connected
components of
$W^u(p)$ for $p\in\cD\cap\Omega$ and of ${\bf R}^2\cap f^{n}\{y=-1\}-I_-$ for
$n\ge0$. 
\end{thm}

\proof  We have seen that $\Omega_+\cap R_j\cap{\bf R}^2$ is homeomorphic to
$\cT_{+,j}\times I$, and thus $\Omega_+\cap{\bf R}^2$ is a lamination of $(R_0\cup
R_1)\cap{\bf R}^2$.  Now $f^{-n-1}\Omega_+\supset f^{-n}\Omega_+$, and
$f^{-n-1}\Omega_+-f^{-n}\Omega_+\subset\interior R_-$.  Thus
$f^{-1}\Omega_+\cap{\bf R}^2$ is a lamination except at the image of the critical
locus of $f^{-1}$, which is $I(f)$.  The backward orbit of $I(f)$ is $I_+$.  Thus
this theorem follows from Theorem 10.5.  
\qed
\begin{figure}
\centerline{\epsfxsize2.6in\epsfbox{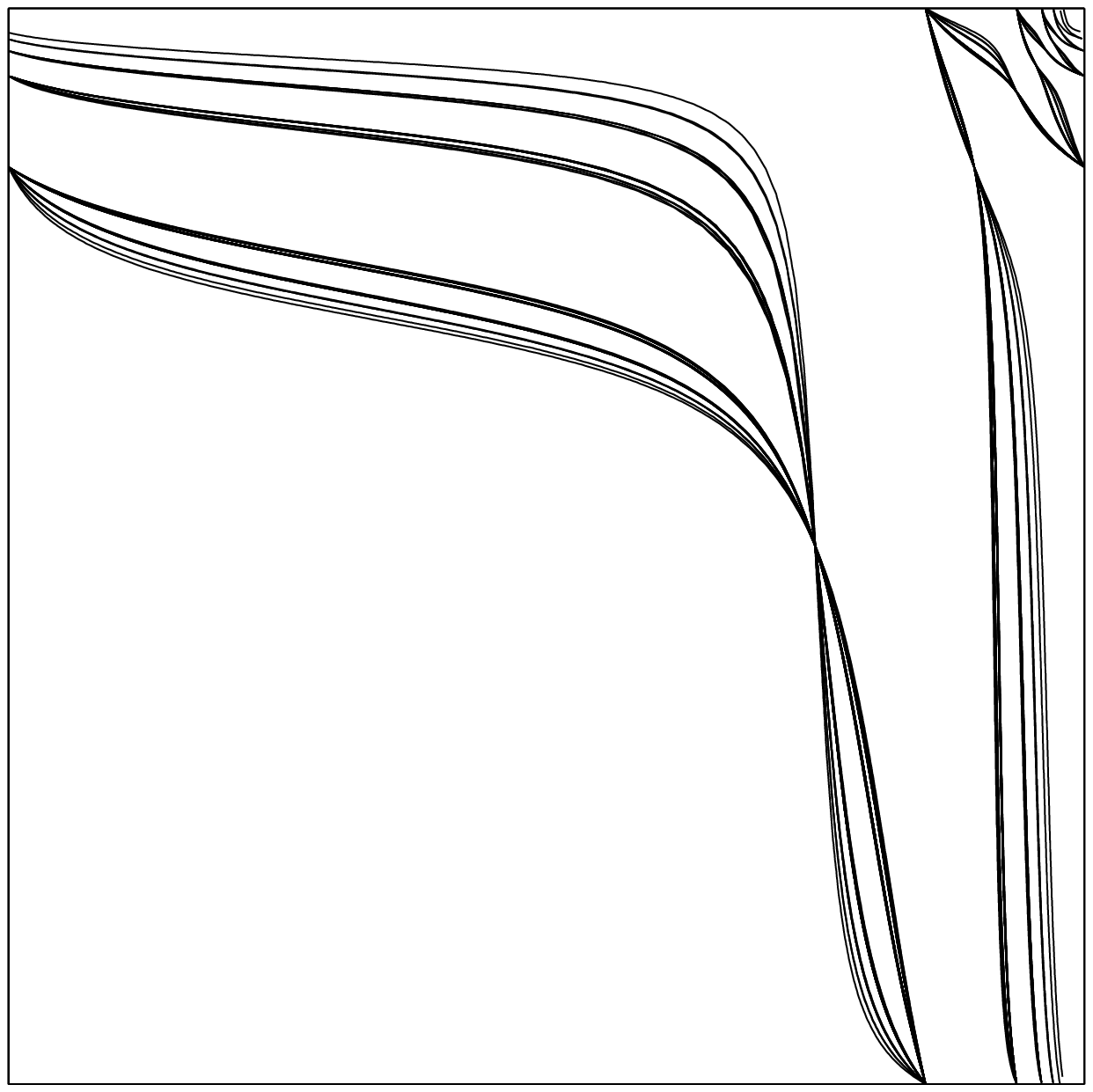}\epsfxsize2.6in
            \epsfbox{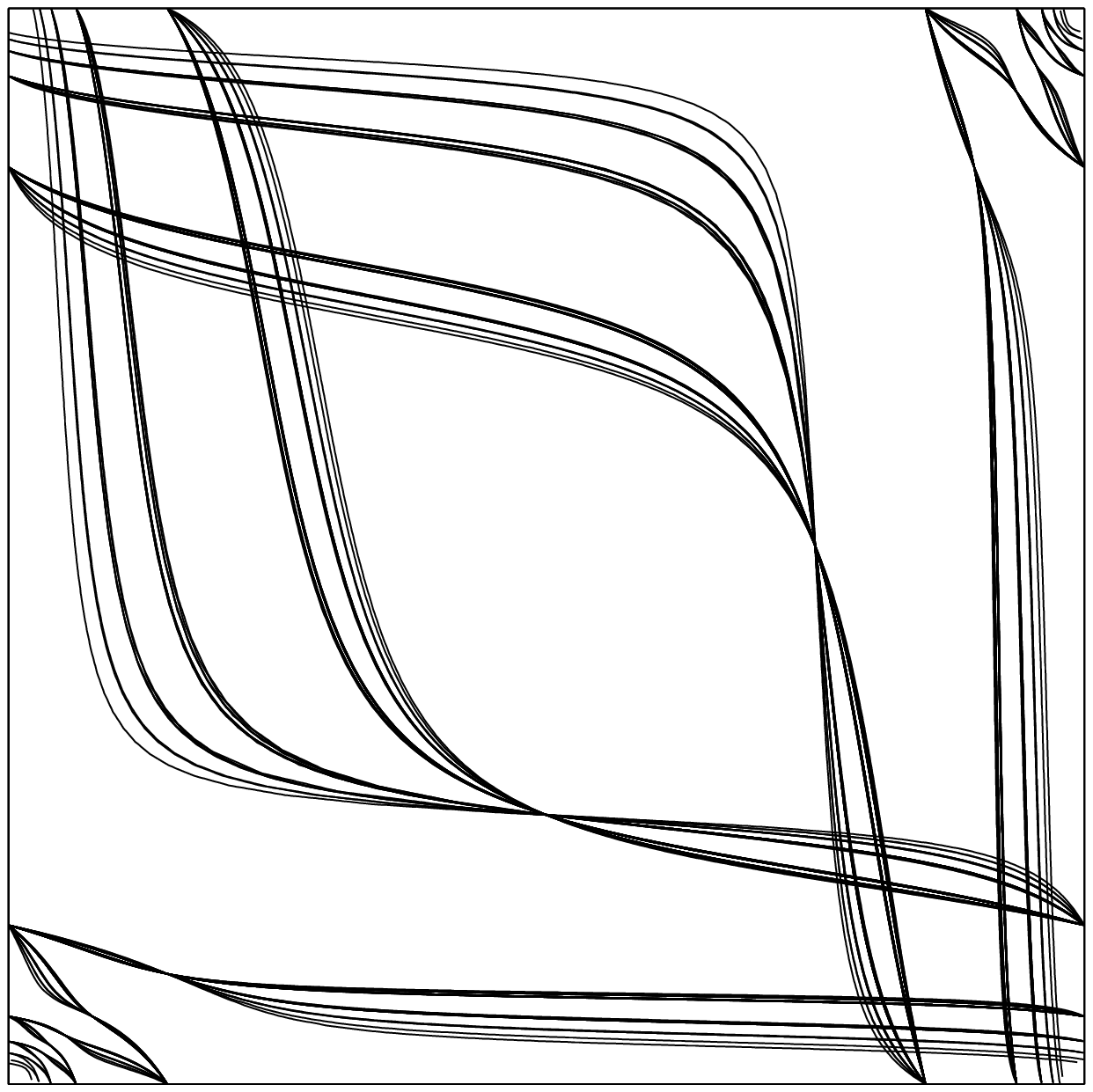}}
\caption{Stable lamination (left), both laminations (right),
         $a=-2$. \label{omegapprox}}
\end{figure}

The left side of Figure \ref{omegapprox}, which shows $f^{-10}(\{x=1.5\})$,
gives an approximation of $\cW^s$.  Note that we have replaced the usual 
Euclidean coordinates $(x,y)$ on $\R^2$ with new coordinates 
$(x',y') = (\arctan x, \arctan y)$ in order to better show the behavior near
infinity.
Note also that the image curve intersects itself precisely at points in 
$I(f^{10})$, and these are in the backward orbit of $\{(-a,\infty), (1,0)\}$.
Both $\cW^s$ and $\cW^u$ are given on 
the right hand side of Figure \ref{omegapprox}.  The intersection 
of the two curves approximates $\Omega$. 

Now we introduce some standard terminology from the theory of currents.  A good
reference for this is Morgan [M].  For an
oriented arc $\gamma$ of finite length, we define the \emph {current of
integration} $[\gamma]$, which is an object that acts on a 1-form $\xi$ as
follows
$$\langle[\gamma],\xi\rangle:=\int_\gamma\xi.$$
Since $\gamma$ has finite length, $\langle[\gamma],\xi\rangle$ is dominated by the
length of $\gamma$ and  $\sup_\gamma|\xi|$.  In this case, $[\gamma]$ is said to be
\emph {represented by integration}.  Currents which are represented by integration
may be treated as vector-valued measures in the sense that they may be written  the
form
$\vec t\cdot\lambda$, where $\vec t$ is a Borel measurable 1-vector, and $\lambda$
is a Borel measure.  Further,  it follows that for any smooth function
$h$,
$\langle [\gamma],dh\rangle=h(Q)-h(P)$, where $P$ and $Q$ are the endpoints of
$\gamma$.  In other words, the boundary of the current $[\gamma]$ is a
difference of point masses:
$\partial[\gamma]=\delta_Q-\delta_P$.
 
It is an elementary consequence of our slope bounds that:
\begin{prop}
The length of $R(w)$ as a curve in 
$\overline{{\bf R}^2}$ is uniformly bounded for all
$w\in\Sigma_+$.
\end{prop}

\begin{figure}
\centerline{\epsfxsize4.6in
            \epsfbox{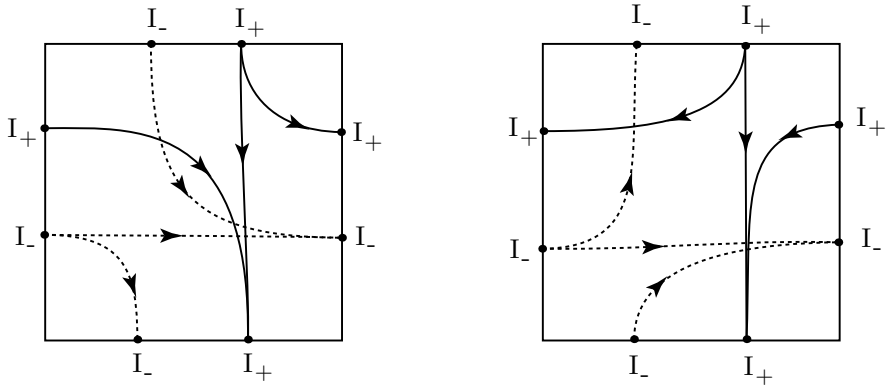}}
\caption{Orientations of stable and unstable currents. \label{orientation}}
\end{figure}

For $w\in\cT^+_j$, we assign to the arc
$R(w)\cap{\bf R}^2$ the orientation for which it  points down and to the right. 
A schematic diagram of this orientation is given for $a<-1$ on the left hand side
of Figure 10.  Note that this diagram in fact shows how to orient the whole
stable lamination
$\cW^s$.  The orientation for $\cW^u$, obtained by mapping under
$(x,y)\mapsto(-y,-x)$, is also given in Figure 10.

The current $[R(w)\cap{\bf R}^2]$ depends continuously on $w$, and we define 
$$\mu_j^+=\int_{w\in\cT_{+,j}}[R(w)\cap{\bf R}^2]\,\nu^+(w).$$
The action of $\mu^+_j$ on a 1-form $\xi$ is given by:
$$\langle \mu^+_j,\xi\rangle:= \int_{w\in\cT_{+,j}} \langle
[R(w)\cap{\bf R}^2],\xi\rangle\,\nu_+(w).$$
In other words, this is a direct integral of currents of integration: for
each
$w\in\cT_{+,j}$, we apply the current of integration $[R(w)]$ to the test form
$\xi$ on ${\bf R}^2$, and we then integrate the result with respect to $\nu_+$ over
the set
$w\in\cT_{+,j}$.

\begin{prop} $\mu^+_j$ is a current on ${\R^2}$ whose support is
contained in 
$\Omega^+\cap R_j$.  The support of its boundary, $\partial\mu^+_j$ is contained
in $\partial^u R_j$.
\end{prop}

Now let us see how $\mu^+_0$ transforms under $f$. The pull-back of the direct
integral is the direct integral of the pull-back.  Therefore we obtain
$$f^*\mu^+_0=\int_{w\in\cT^+_0} 
f^*[R(w)\cap{\bf R}^2]\,\nu_+(w) =\int_{w\in\cT^+_0} 
[f^{-1}(R(w)\cap{\bf R}^2)]\,\nu_+(w).$$
We compare $f^*\mu^+_j$ with $\mu^+_j$.  Since $\mu^+_j$ is represented by
integration, we may multiply it by ${\bf 1}_{R_0\cup R_1}$, the characteristic
function of the set $R_0\cup R_1$.   We know that for each
$w$,
$f^{-1}(R(w)\cap{\bf R}^2)$ crosses both
$R_0$ and
$R_1$ in s-arcs, as is shown in Figure 8.  Figure 8 shows that the orientations
of
$f^{-1}(R(w)\cap {\bf R}^2)$ are opposite from the orientations in Figure 10.  By
the definition of $\sigma^+_k$, we see that 
$${\bf 1}_{R_0\cup R_1} [f^{-1}(R(w)\cap{\bf R}^2)] =
-[R(\sigma^+_0w)]-[R(\sigma^+_1w)].$$
If we take the integral with respect to $\nu_+$, we have
$${\bf 1}_{R_0\cup
R_1}f^*\mu^+_0 = -\int_{w\in\cT_{+,0}}\nu_+(w)\,([R(\sigma^+_0w)] +
[R(\sigma^+_1w)]).$$
By the balanced property of $\nu_+$, this is
$$= -\phi\int_{v_0\in\sigma^+_0(\cT_{+,0})\subset\cT_{+,0}} \nu_+(v_0)\,[R(v_0)]
-\phi\int_{v_1\in\sigma^+_1(\cT_{+,0})=\cT_{+,1}} \nu_+(v_1)\, [R(v_1)].$$
Applying the same reasoning to ${\bf 1}_{R_0\cup R_1}f^*\mu_1^+$, and adding, we
obtain:
\begin{prop}
$${\bf 1}_{R_0\cup R_1}f^*(\mu^+_0+\mu^+_1) = -\phi (\mu^+_0+\mu^+_1).$$
\end{prop}
Now we consider the sequence of currents $(-\phi)^{-n}f^{*n}(\mu^+_0+\mu^+_1)$,
which coincide with $\mu^+_0+\mu^+_1$ on the interior of $R_0\cup R_1$.  We may
now argue as in the proof of Theorem 10.8 to obtain:
\begin{thm}
The current $\mu^+_0+\mu^+_1$ extends to a closed current $\mu^+_{\bf R}$ on
${{\bf R}^2}-I_+$, which satisfies $f^*\mu^+_{\bf
R}=-\phi\cdot\mu^+_{\bf R}$.
\end{thm}

In a similar way, we may define currents
$$\mu^-_j=\int_{w\in\cT_{-,j}}[R(w)\cap{\bf R}^2]\,\nu^-(w),$$
and we may extend $\mu^-_0+\mu^-_1$ to a closed current $\mu^-_{\bf R}$ on ${\bf
R}^2-I_-$, which satisfies $f^*\mu^-_{\bf R}=-\phi\mu^-_{\bf R}$.

Let $\gamma_1$ and $\gamma_2$ be oriented arcs of finite length, and let
$[\gamma_1]$ and $[\gamma_2]$ be their currents of integration.  We define the
wedge (intersection) product of these currents to be
$[\gamma_1]\wedge[\gamma_2]=0$ if the supports of
$[\gamma_1]$ and $[\gamma_2]$ are disjoint.  If $\gamma_1\cap\gamma_2$ consists of
a single point $p$, where the curves intersect transversally, we define
$[\gamma_1]\wedge[\gamma_2]=\pm\delta_p$ to be the (signed) point mass at $p$, with
the choice of sign taken to be positive if wedge products of the tangents at $p$
agrees with the orientation of ${\bf R}^2$ at $p$.  Let us orient $R_j$ such
that for $w\in\Sigma'$ we have
$[R(w^+)]\wedge [R(w^-)]=\delta_{R(w)}.$

The wedge product of the direct integrals is the direct integral of the wedge
products, so we obtain
\begin{eqnarray}
 \mu^+_{\R}\wedge\mu^-_{\R} =  &
                \int_{w\in\Sigma_+}\int_{\tilde w\in\Sigma_-} [R(w)]\wedge
                                                           [R(\tilde w)] 
                        \,\nu^+(w)\otimes\nu^-(\tilde w)\nonumber\\
		= &\int_{\tilde w_0 = w_0} \delta_{R(w)\cap R(\tilde w)}
                           \,\nu^+(w)\otimes\nu^-(\tilde w)\ \ \ \ \ \ \ \ \ \ \ \ 
\nonumber
\end{eqnarray}             

From \S2 we recall the product structure map
$$\cT_{+,j}\times\cT_{-,j}\ni(w^+,w^-)\mapsto\tilde w\in\Sigma'\cap\pi^{-1}(j),$$
defined by $\tilde w^\pm=w^\pm$.  With this map, we have the relation
$\nu^+\otimes\nu^-\cong\nu$.  In \S8 we discussed the product structure mapping
$$R:\cT_{j,+}\times\cT_{j,-}\ni(w^+,w^-)\mapsto R(\tilde w)\in\Omega\cap
R_j\cap{\bf R}^2.$$
This map was used to define $\mu$ as the pushforward of $\nu=\nu^+\otimes\nu^-$. 
Since these two product structures coincide, we have:
\begin{thm}
$\mu = \mu^+_{\R}\wedge\mu^-_{\R}$.
\end{thm}

\section{Parameter values $-1<a<0$}
\label{smallnega}
\begin{figure}
\centerline{\epsfxsize2.6in
\epsfbox{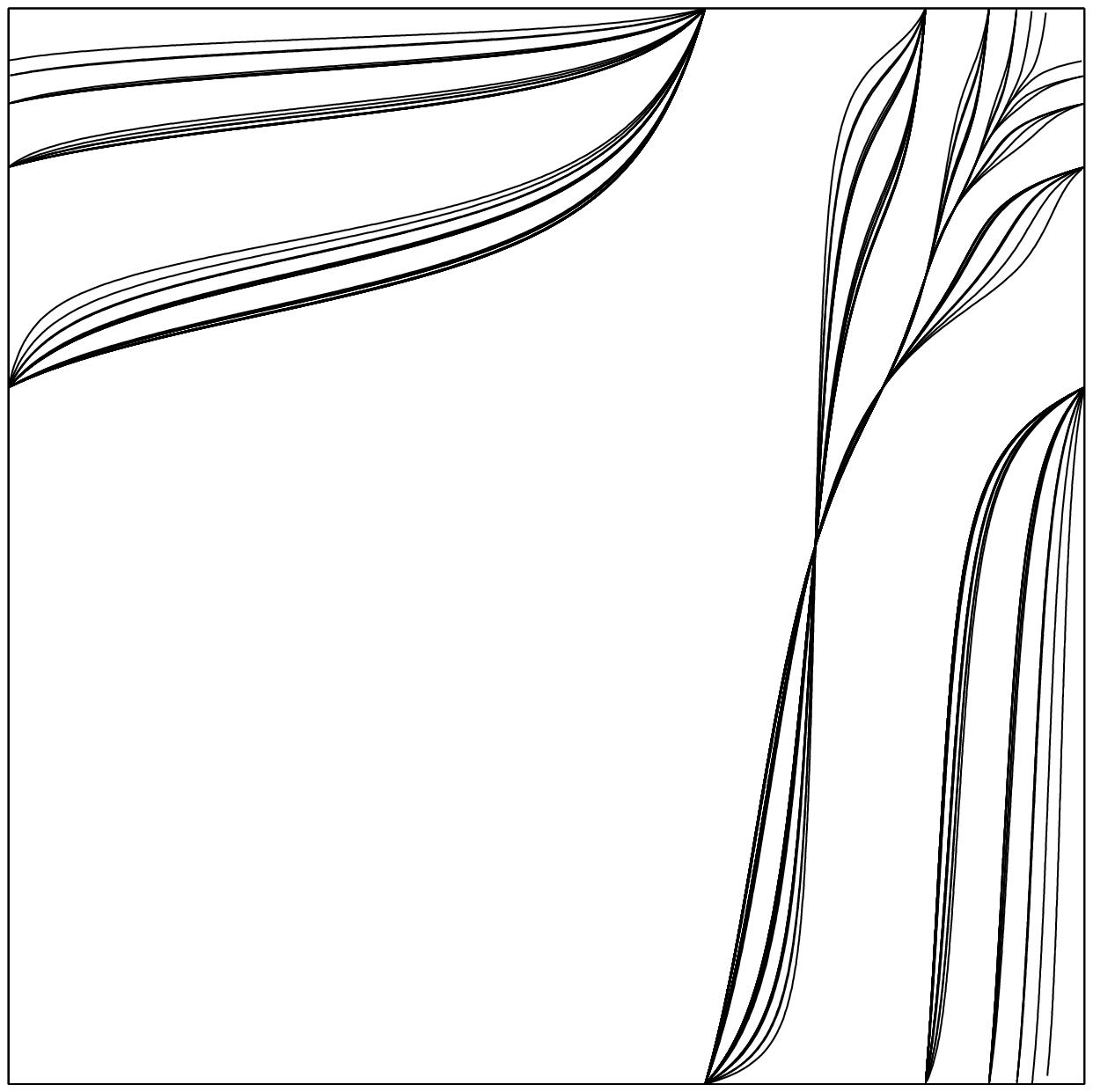} 
\epsfxsize2.6in\epsfbox{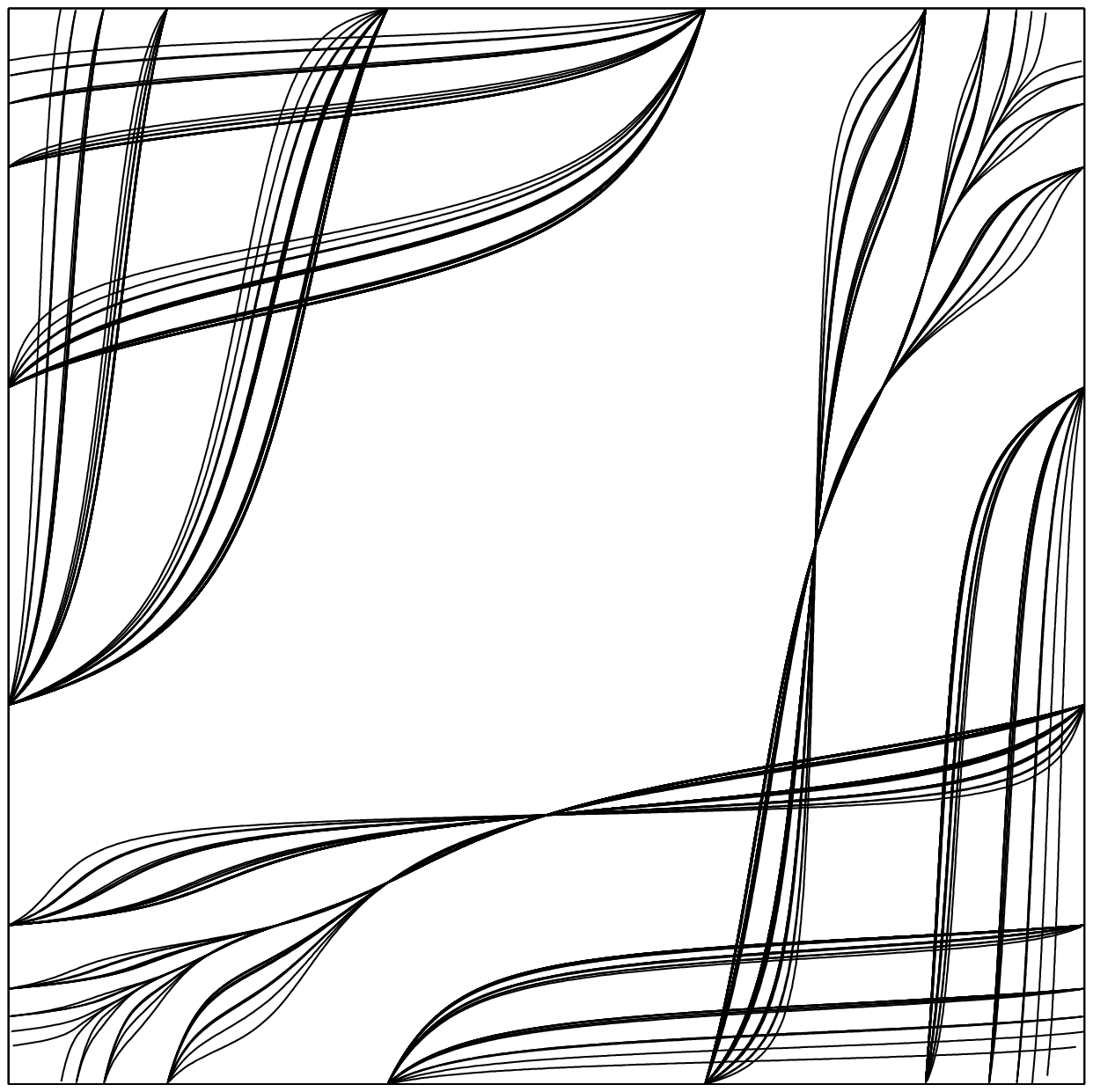}}
\caption{Stable lamination (left), both laminations (right), $a=-.5$
\label{omegapprox2}}
\end{figure}

\noindent The main results that we have proved so far hold for all parameters 
$a < 0$ save the exceptional value $a=-1$.  However, in our proofs we have 
been assuming that $a<-1$.  The details in the case $-1<a<0$ are similar
enough that most of them are not worth repeating.  We will use this section
to point out those few places in which the differences are significant.
\begin{figure}
\centerline{\epsfxsize3.6in
\epsfbox{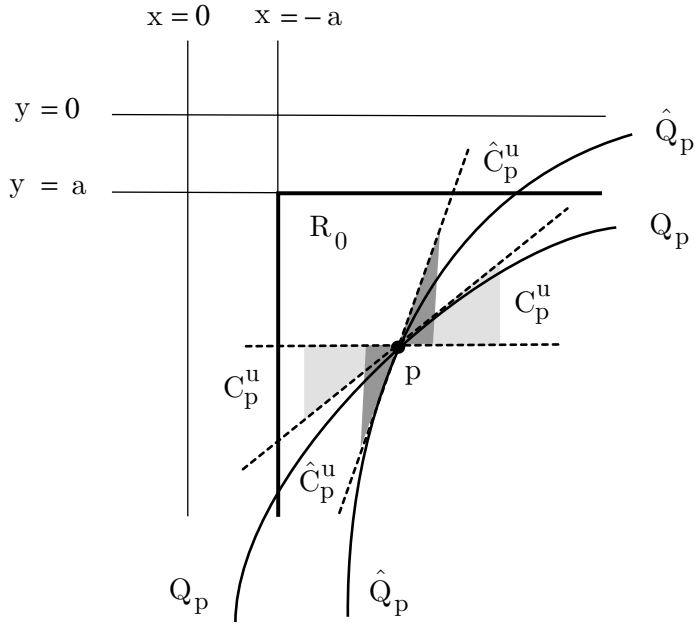} }
\caption{Cone Fields for $-1<a<0$
\label{cones2}}
\end{figure}

The first and most basic change is that the definition of the filtration
$R_0,R_1,R_\pm$ must be changed so that the lines $\{x=-a\}$ and $\{y=a\}$
take the place of the lines $\{x=1\}$ and $\{y=-1\}$.  In other words, $R_1$ is
unchanged, but $R_0=\{-a\le x\}\cap\{y\le a\}$.  Proposition \ref{thm1} then holds
with a slight alteration to the first two items.  For example, the first item 
should read: \emph{if $(x_0,y_0)\in R^+$, then}
$$
\max\{x_1 + a, y_1\} \leq \max\{x_0+a,y_0\} + a.
$$
The definitions of s-arcs and u-arcs do not change.  
Nor does Proposition \ref{allowables}.  
We define the stable/unstable cones as follows.  For $p\in R_0$, we let $Q_p$
denote the branch of the hyperbola (quadric) passing through $p$ and having
$\{x=0\}$ and $\{y=a\}$ as asymptotes.  We let $\hat Q_p$ denote the hyperbola
passing through $p$ and having $\{x=a\}$ and $\{y=0\}$ as asymptotes.  For $p\in
R_0$, we define
$\cC^u_p$ (respectively, $\hat \cC^u_p$) to be the cone of vectors swept out by
starting at the horizontal $H_p$ and moving counter-clockwise to the tangent to
$Q_p$ (respectively $\hat Q_p$) at $p$.  This is shown in Figure \ref{cones2}. 
The following is proved along the same general lines as Theorems 5.1 and 5.2.
\begin{thm}
If $p,fp,f^2p\in (R_0\cup R_1)\cap {\bf R}^2$, then
$Df_p\hat\cC^u_p\subset\cC^u_{fp}$.  Further, $Df^2_p$ maps $\cC^u_p$ strictly
inside
$\cC^u_{f^2p}$.
\end{thm}
From this theorem, we derive the  slope bounds analogous to those
of Theorem \ref{slope}. 
We change the definitions of {uniform} s-arcs and u-arcs to fit the new slope
bounds.  Despite the change in particulars,
the role of the stable and unstable cones is exactly as before.  For instance,
the angle of intersection between a uniform $s$-arc and a uniform $u$-arc is
bounded away from zero uniformly on any compact subset of $\R^2$.

All of the assertions in Sections  6--10 remain true without change for $-1<a<0$. 
However, we note that in Section \ref{parabolic} some of the proofs become a little
easier.  The reason for this can be seen in the approximation of $f^2$ 
presented in Proposition \ref{approx}.  For $x$ and $y$ large and $a$
between $0$ and $-1$, we have that $x$ decreases and that $y$ increases
with every iterate of $f^2$.  This was not quite true in the case $a<-1$.

The most important difference between the cases $a<-1$ and $-1<a<0$ concerns
orientation.  We orient uniform $u$-arcs and $s$-arcs so that their
tangent vectors point up and to the right; a schematic illustration is given on
the right hand side of Figure 10.  In this case, $f$ preserves, rather than
reverses, their orientations.  Thus we have the transformation law 
$$
f^*\mu_{\R}^+|_{R_0\cup R_1} = \phi\mu^+_{\R}.
$$ 
for the current $\mu^+_{\R}$ supported on $\Omega^+$.

\section{The purely complex point of view}
\label{complex}
\noindent We have given a detailed description of the dynamics of the specific family
$\{f_a:a<0,a\ne-1\}$ of birational maps of ${\bf R}^2$.  There is also a general theory which
applies to birational mappings of ${\bf C}^2$.  In this Section we state a result of [DF] on
the existence of invariant currents and two results from [BD] on the existence of invariant
measures.  Then we show how the results obtained in the preceding sections fit into the
framework of the complex theory.

From \S3 recall the basis $\{\gamma_1^*,\gamma_2^*\}$ of $H^2({\bf P}^1\times{\bf P}^1)$ and
the action of $f^*$ on $H^2({\bf P}^1\times{\bf P}^1)$.  The $\phi$ eigenspace of $f^*$ is
generated by the cohomology class $\theta^+:=c(\gamma^*_1+\phi\gamma^*_2)$ for any $c\ne0$. 
The generator $\theta^-$ for the $\phi$ eigenspace of $f_*=(f^{-1})^*$ is given by
$c(\phi\gamma^*_1+\gamma^*_2)$.  We choose $c=(1+\phi^2)^{-1/2}>0$ so that
$\theta^+\cdot\theta^-=1$.  Because
$f^*$ and $f_*$ are adjoint with respect to the intersection product, standard linear algebra,
applied to the 2-dimensional space $H^2$, gives
$$
\lim_{n\to\infty}\frac{f^{n*}\theta}{\phi^n} = (\theta\cdot\theta^-)\theta^+
$$
for any class $\theta\in H^{2}(\pp)$.  

A current $T$ of bidimension (1,1) is said to be {\em positive} if
$\pair{T}{\frac1i\alpha\wedge\bar\alpha}
\geq 0$  for any $(1,0)$ form $\alpha$ on $\pp$.  $f$ induces a well-defined action, also
denoted by $f^*$, on the space of positive, closed currents (see [Gu1] or [S]).  In fact, if
$\omega$ is a current representing a cohomology class
$\theta$, then
$f^*\omega$ is a current representing the cohomology class $f^*\theta$.  If 
$\omega$ is smooth, then $f^{n*}\omega$ is smooth except at the finite set $I(f^n)$.  The
action of $f^*$ on positive, closed currents closely follows the action of $f^*$ on cohomology:
\begin{thm}
There exists a unique positive, closed $(1,1)$ current $\mu^+$ on $\pp$ with the 
following properties.  
\begin{itemize}
\item $\mu^+$ represents $\theta^+$;
\item $f^*\mu^+ = \phi\mu^+$;
\item For every cohomology class $\theta\in H^2(X)$ and every smooth form
      $\omega$ representing $\theta$,
$$
\lim_{n\to\infty} \frac{f^{n*}\omega}{\phi^n} = (\theta\cdot\theta^-)\mu^+.
$$
\end{itemize}
\end{thm}

Applying the previous theorem to $f^{-1}$, we obtain an invariant current $\mu^-$ in the
cohomology class of $\theta^-$.  Thus, $\mu^-=\phi^{-1}(f^{-1})^*\mu^-=\phi^{-1}f_*\mu^-$.  
Because 
$$
\overline{I_+}\cap \overline{I_-}
= \{(\infty,\infty)\}
$$ 
contains only one point, we also have 
\begin{thm}
The wedge product $\mu = \mu^+\wedge\mu^-$ is well-defined and equal to an 
$f$-invariant ergodic probability measure on $\pp$.
\end{thm}

It is useful to know that  
the limit in the third item of the first theorem can be freely interchanged 
with the wedge product defining $\mu$.  
Thus the measure $\mu$ arises from simultaneously pushing forward and pulling
back arbitrary smooth currents.
\begin{thm}
\label{cpxconvergence}
Let $\theta,\tilde\theta\in H^{2}(\pp)$ be classes represented by smooth currents
$\omega,\tilde\omega$.  Then
$$
c\mu = \lim_{n,m\to\infty} 
      \frac{f^{n*}\omega\wedge f^m_*\tilde\omega}{\phi^{n+m}},
$$
where $c = (\theta\cdot\theta^-)(\tilde\theta\cdot\theta^+)$.
\end{thm}

Let us consider the cohomology class $\gamma_2^*$ generated by any vertical $\{x=t\}$.  If
$\chi\in C_0^\infty({\bf R})$, $\int\chi=c$, then
$$\omega:=\int[V_t]\,\chi(t)\,dt$$
is a current which represents $c\gamma_2^*$.  Similarly, if $\tilde\chi\in C_0^\infty({\bf
R})$, $\int\tilde\chi=\tilde c$, and if
$H_s=\{y=x\}$, then
$$\tilde\omega:=\int[H_s]\,\tilde\chi(s)\,ds$$
is a  current representing the cohomology class $\tilde c\gamma_1^*$ of a horizontal line.  We
choose $c,\tilde c>0$ so that
$(\theta^-\cdot\omega)(\theta^+\cdot\tilde\omega)=1$.  While $\omega$ and $\tilde\omega$ are
not smooth on ${\bf P}^1\times{\bf P}^1$, they have continuous potentials, and Theorem 12.3 can
be adapted to cover this case.  It follows, then, that
$\phi^{-(n+m)}f^{*n}\omega\wedge f^m_*\tilde\omega$ converges to the measure $\mu$ of Theorem
12.2 as $n,m\to\infty$.   

\begin{thm}
 Suppose that $a<0$, $a\ne -1$.  Then for any $s\le0$, $t\ge1$, we have
$$\mu=\frac{1+\phi^2}{\phi^2}\lim_{n,m\to\infty}\phi^{-n-m}\sum_{a\in f^{-n}V_s\cap
f^mH_t}\delta_a,$$ where the convergence is in the weak sense of measures.
\end{thm}

\proof  
With the notation as above, let us suppose that $\chi$ is
supported in $\{1\le x<\infty\}$, and thus for each $t$ in the integral defining
$\omega$, $V_t$ intersects
$R_0\cup R_1$ in a uniform s-arc.  Similarly, if $\tilde\chi$ is supported in
$\{-\infty<y\le-1\}$, then for each $s$ in the integral defining $\tilde\omega$, 
$H_s$ intersects $R_0\cup R_1$ in a uniform u-arc.  Now let us consider fixed $n,m>0$, and let
$w$ be a word of extent $[-n,m]$.  By
\S4, it follows that
$f^{-n}V_t\cap f^mH_s$ intersects $R(w)$ in exactly one point, and the union over all $w$ gives 
$f^{-n}V_t\cap f^mH_s$.
Let us denote this point by $p(s,t,w)$.  It follows that the restriction of
$f^{m*}\omega\wedge f^n_*\tilde\omega$ to $R(w)$ is given by the integral
$$I_w:=\int_{s\in{\bf R}}\int_{t\in{\bf R}}\delta_{p(s,t,w)}\,\chi(s)\tilde\chi(t)dsdt.$$

Next let us note that as a consequence of Theorem 8.1, we have the following:  if $K$ is a
compact subset of ${\bf R}^2$, then
$$\lim_{m,n\to\infty}\max\{diam(R(w)\cap K)\}=0,$$
where the maximum is taken over all words $w$ of extent $[-n,m]$.  By Theorem 12.3, we know
that the sum of the integrals $\sum_wI_w$ converges to $\frac{\phi^2}{1+\phi^2}\mu$ as
$n,m\to\infty$.  Since the diameters of the $R(w)$ shrink to zero uniformly on any compact set
$K$, it follows that the difference between 
$\sum_wI_w$ and $\phi^{-n-m}\sum_{a\in f^{-n}V_s\cap f^mH_t}\delta_a$
tends to 0 as $n,m\to\infty$.  This proves the Theorem.
\qed

As was seen in the proof of the preceding theorem, each measure $I_w$ has the same mass, and so 
$f^{n*}\omega\wedge f^m_*\tilde\omega$ puts the same mass on each
the rectangle $R(w)$.  On the other hand, a well known property of the measure $\nu$ is that
it can be obtained by equidistributing mass over the ``cylinder sets''
$$
C(w) = \{ s\in\Sigma:s_j = w_j \text{ for } -n\leq j\leq m\}
$$ 
and letting $n,m\to\infty$.  Now let us recall that the measure $\mu$ in \S8 was obtained as
the image of $\nu$ under the map $R$.  Since $R$ takes a cylinder $C(w)$ to a rectangle
$R(w)$, we have the following:

\begin{thm}
\label{samemeasure}
Suppose that $a<0$, $a\neq -1$.  Then the measure $\mu$ defined in this section 
is the same as the measure $\mu$ defined in \S\ref{conjugacy}.
\end{thm}

\bigskip
\centerline{\bf References}

{[Ab1]}  N. Abarenkova, J.-Ch. Angl\`es d'Auriac, S. Boukraa,
S. Hassani and J.-M. Maillard, Rational dynamical zeta functions for
birational transformations, Physica A 264 (1999) 264--293.

{[Ab2]}  N. Abarenkova, J.-Ch. Angl\`es d'Auriac, S. Boukraa,
S. Hassani and  J.-M. Maillard, Topological entropy and complexity for
discrete dynamical systems, Phys. Lett. A 262 (1999) 44--49.

{[Ab3]}  N. Abarenkova, J.-Ch. Angl\`es d'Auriac, S. Boukraa,
S. Hassani and 
J.-M. Maillard, Growth complexity spectrum of some discrete dynamical
systems, Physica D 130 (1999) 27--42.

{[Ab4]}  N. Abarenkova, J.-Ch. Angl\`es d'Auriac, S. Boukraa  and 
J.-M. Maillard, Real Arnold complexity versus real
topological entropy for birational transformations, J. Phys. A. 33
(2000), 1465--1501.

{[Ab5]}  N. Abarenkova, J.-Ch. Angl\`es d'Auriac, S. Boukraa and 
J.-M. Maillard, Real topological entropy versus metric entropy for
birational measure-preserving transformations, Physica D 144 (2000)
387--433.

[Ak] E. Akin, {\sl The General Topology of Dynamical Systems}, AMS,
1993.

[BD] E. Bedford and J. Diller, Energy and invariant measure for bimeromorphic surface maps.

{[BLS]} E. Bedford, M. Lyubich, and J. Smillie, Polynomial diffeomorphisms
of ${\bf C}^2$, IV. The measure of maximal entropy and laminar currents. Invent.
Math., 112 (1993), 77--125.

[BS1]  E. Bedford and J. Smillie, Real polynomial diffeomorphisms with maximal entropy: Tangencies.
http://arXiv.org/math.DS/0103038

[BS2] E. Bedford and J. Smillie, Real polynomial diffeomorphisms with maximal entropy: II. Small Jacobian.

[BTR] M. Bernardo, T.T. Truong and G. Rollet, The discrete Painlev\'e I equations:
transcendental integrability and asymptotic solutions, J. Phys. A: Math. Gen. 34 (2001),
3215--3252.

[BM] S. Boukraa and J.-M. Maillard, Factorization properties of birational mappings, Physica
A 220 (1995), 403--470.

[dCH]  A. de Carvalho and T. Hall, How to prune a horseshoe, Nonlinearity, 15 (2002), no. 3,
R19--R68. 

[DN]  R. Devaney and Z. Nitecki, Shift automorphisms in the H\'enon family, Comm. Math. Phys. 67 (1979),
137--146.

[DF]   J. Diller and C. Favre,  Dynamics of bimeromorphic
maps of surfaces, Amer. J. Math., 123 (2001), 1135--1169.

[Du] R. Dujardin,  Dynamique d'applications non polynomiales et courants laminaires, Doctoral
Thesis, Univ. de Paris Sud, Orsay, 2002.

[Fr]  S. Friedland, Entropy of rational self-maps of projective
varieties. In {\sl Dynamical Systems and Related Topics} (Nagoya, 1990),
pages  128--140. World Sci.\ Publishing, River Edge, NJ, 1991.

[FM] S. Friedland and J. Milnor, Dynamical properties of plane polynomial automorphisms, Ergod.\
Th.\ \&\ Dynam.\ Sys., 9 (1989),  67--99.

[Fu]  W. Fulton, {\sl Intersection Theory}, Ergebnisse der Mathematik und ihrer Grenzgebiete,
Springer-Verlag, 1984.

[GNR] B. Grammaticos, F. Nijhoff and A. Ramani, Discrete Painlev\'e equations.  The
Painlev\'e property.  (CRM Series in Mathematics Physics), (New York: Springer) (1999),
413--516.

[Gu1] V. Guedj, Dynamics of polynomial mappings of ${\bf C}^2$, Am. J. Math., 124 (2002),
75--106.

[Gu2] V. Guedj,  Rational mappings with large topological degree, preprint.

[HO] J.H. Hubbard and R. Oberste-Vorth, H\'enon mappings in the complex domain II: Projective
and inductive limits of polynomials, {\sl Real and Complex Dynamical Systems} (B. Branner and
P. Hjorth, eds.), Kluwer, Boston, 1995, pp.\ 89--132.

[KH] A. Katok and B. Hasselblatt, {\sl Introduction to the Modern Theory of Dynamical Systems},
Cambridge U. Press (1995).

[LM] D. Lind and B. Marcus, {\sl Symbolic Dynamics and Coding}, Cambridge U. Press (1995).

[M]  F. Morgan, {\sl Geometric Measure Theory}, Academic Press, 1988.

[RS]  D. Ruelle and D. Sullivan, Currents, flows and diffeomorphisms,
Topology, 14 (1975), 319--327.

[S] N. Sibony, Dynamique des applications rationnelles de ${\bf P}^k$, Dynamique et
g\'eom\'etrie complexes (Lyon, 1997), 97--185,  Panor. Synth\`eses, 8, Soc.\ Math.\ France,
Paris, 1999.

\nocite{GrHa}

\bibliographystyle{texfiles/mjo}
\bibliographystyle{plain}
\bibliography{texfiles/references}

\end{document}